\documentclass{article}
\usepackage{amsmath,amsthm,amsfonts,amscd,amssymb,eucal,latexsym}

\begin{document}
\baselineskip15pt
\textwidth=12truecm
\textheight=20truecm
\hoffset-.1cm
\voffset-.5cm

\def\pp{\psi\vphantom A}
\def\ff{\varphi\vphantom A}
\def\ov{\overline}
\def\O{{\mathcal O}}
\def\var{\varepsilon}
\def\h#1{\hbox{\rm #1}}
\font\got=eufm10 scaled\magstephalf
\def\Ag{\h{\got A}}
\def\Fg{\h{\got F}}
\def\lam{\lambda}
\def\T{{\cal T}}
\def\C{{\cal C}}
\def\A{{\cal A}}
\def\R{{\cal R}}
\def\B{{\cal B}}
\def\M{{\cal M}}
\def\D{{\cal D}}
\def\N{{\cal N}}

\def\Bbb{\mathbb}
\def\R{{\cal R}}
\def\B{{\cal B}}
\def\P{{\mathbb P}}
\def\T{{\mathbb T}}
\def\N{{\cal N}}
\def\C{{\cal C}}
\def\d{\hbox{\rm $\,$d}}
\def\End{\hbox{\rm End}}
\def\Aut{\hbox{\rm Aut}}
\def\O{{\cal O}}
\def\H{{\cal H}}
\def\I{{\cal I}}
\def\cC{{\cal C}}
\def\E{{\cal E}}
\def\cT{{\cal T}}
\def\N{{\cal N}}
\def\M{{\cal M}}
\def\D{{\cal D}}
\def\cB{{\cal B}}
\def\K{{\cal K}}
\def\A{{\cal A}}
\def\U{{\cal U}}
\def\Z{{\mathbb Z}}
\def\ov{\overline}
\def\Om{\Omega}
\def\var{\varepsilon}
\font\got=eufm10 scaled\magstephalf
\def\Ag{\h{\got A}}
\def\Fg{\h{\got F}}
\def\lam{\lambda}

\font\eightrm=cmr8

\title{An Algebraic Duality Theory for Multiplicative Unitaries}
\author{Sergio Doplicher, Claudia Pinzari\\
Dipartimento di Matematica, Universit\`a di Roma ``La Sapienza,''\\ 
I-00185 Rome, Italy\\
\\
 and John E.\ Roberts\\ 
Dipartimento di Matematica, Universit\`a di Roma ``Tor Vergata,''\\
I-00133 Rome, Italy\\
}
\date{}
\maketitle

\begin{abstract}
Multiplicative unitaries are described in terms of a pair of 
commuting shifts of relative depth two. They can be generated 
from ambidextrous Hilbert spaces in a tensor $C^*$--category.
 The algebraic analogue of the Takesaki-Tatsuuma Duality Theorem 
characterizes abstractly $C^*$--algebras acted on by unital endomorphisms
that are intrinsically related to the regular representation of a 
multiplicative unitary. The relevant $C^*$--algebras turn out to 
be simple and indeed separable if the corresponding multiplicative 
unitaries act on a separable Hilbert space.
A categorical analogue provides internal characterizations of minimal
representation categories of a multiplicative unitary. Endomorphisms of
the Cuntz algebra related algebraically to the grading are discussed as 
is the notion of braided symmetry in a tensor $C^*$-category.

\end{abstract}

\vfill

\noindent Research supported by MURST and CNR--GNAFA.\eject

\section{Introduction}
{}
An interesting open problem in the duality theory of tensor 
$C^*$--categories is to decide which ones admit a (tensor--preserving) 
embedding into the tensor $C^*$-category of Hilbert spaces. A first 
positive result in this direction is the duality theorem of \cite{DR} asserting 
that symmetric tensor $C^*$--categories with conjugates and irreducible 
tensor unit admit such an embedding.  On the other hand, the theory 
of dimension introduced in \cite{LR} allows one to see that certain tensor 
$C^*$--categories with conjugates and irreducible tensor unit cannot  
be embedded. For such tensor $C^*$--categories admit an intrinsic 
dimension function defined on objects. Under certain circumstances, 
the embedding functor must preserve dimensions. This is well known 
to be so in the rational case, i.e.\ when the set of equivalence classes of 
irreducibles is finite, since it is a simple consequence of the Perron--Frobenius 
Theorem. However, it is also true in the amenable case\cite{LR}. We are 
in the amenable case whenever the category admits a unitary braiding, 
see Theorem 5.31 of \cite{LR}. This means that the tensor $C^*$--categories with conjugates 
and a braiding appearing in low dimensional quantum field theory 
cannot be embedded whenever the dimensions are non--integral. 
On the other hand, the finite dimensional unitary representation theory 
of compact quantum groups provides us with examples of tensor 
$C^*$--categories with conjugates and irreducible tensor units with 
non--integral dimensions which are embeddable by construction. 
The duality theorem of Woronowicz \cite{TK} generalizing Tannaka--Krein to
compact quantum groups provides a way of recognizing such categories if
the embedding is given.\smallskip

   We will not treat the embedding problem here in full generality;
we shall instead present a positive solution for abstract tensor
$C^*$--categories  with specific additional structure. Such structure is
present in a class of model tensor $C^*$--categories namely the minimal
tensor $C^*$--categories generated by the regular representation 
of locally compact groups or of multiplicative unitaries in the sense
of Baaj and Skandalis \cite{BS}. We prove that abstract tensor
$C^*$--categories with this additional structure  are isomorphic to a
model tensor $C^*$--category and are hence embeddable.\smallskip

Our result, Thm.~6.13, is thus a duality theorem for multiplicative
unitaries and is hence applicable to the case of locally compact quantum
groups \cite{KV}. Multiplicative unitaries have already played a role
in Tatsuuma's duality theorem for locally compact groups \cite{T} where
the group elements are identified  in the regular representation using 
the multiplicative unitary and in Takesaki's Hopf von Neumann algebra
version of duality \cite{Tak}. Multiplicative unitaries express the
fundamental property of the regular representation $V$, namely that
$V\times V$ is equivalent to a multiple of $V$. This can also be expressed
through the existence of a remarkable Hilbert space $H$ of intertwiners
from $V$ to $V\times V$ (cf also \cite{D}). Our minimal model is the
smallest tensor
$C^*$--category containing the regular representation as an object and $H$
as a subspace of arrows.\smallskip

A closely related result, Thm.~6.11, characterizes a class of
$C^*$--algebras $\A$ acted on by an endomorphism $\rho$ encoding the
regular representation. Intertwiners between powers of this endomorphism 
encode intertwiners between tensor powers of the regular representation.
The minimal model here is obtained as follows: 
take the  regular representation $V$ considered
 as an object in the 
tensor $C^*$--category
of representations of the multiplicative unitary
and associate to it as in
\cite{DR} the $C^*$--algebra $\O_V$ and its endomorphism
$\rho_V$.   Then the minimal model is the smallest 
$\rho_V$--stable $C^*$--subalgebra containing $H$  acted on by the
restriction of $\rho_V$. The $C^*$--algebra obtained in this way 
is simple and is separable  if and only if the given multiplicative
unitary acts on a separable Hilbert space. 
Further variants of these duality results 
are Theorems 6.5, 6.6, 6.7 and 6.12.\smallskip

In connection with the above results, attention should be drawn to 
Longo's characterization of actions of finite dimensional Hopf algebras
\cite{L} which has a similar algebraic and categorical flavour.
However, his axiomatic structure involving $Q$--systems 
is quite distinct from those used here. Despite this, the Hilbert 
space $H$ puts in an appearance here too and a multiplicative unitary
appears in
his proof.\smallskip

The principal results in the remaining sections may be summarized as
follows. Section 2 gathers together 
some elementary results on categories of Hilbert spaces. The discussion 
centres round the concept of {\it shift}, a category of Hilbert spaces 
with objects labelled by the integers, $0$ being irreducible, and equipped 
with a normal $^*$--functor adding one on objects. Such a structure is 
isomorphic to the category of tensor powers of a given Hilbert space $K$ 
with the functor being tensoring on the right by $1_K$.\smallskip 

  In Section 3, it is shown how two commuting shifts of relative depth 
two are equivalent to giving a multiplicative unitary for the tensor 
product defined by one of the commuting shifts. Duality for multiplicative 
unitaries reflects the symmetry between the commuting shifts. The 
representation category for a multiplicative unitary finds a natural 
expression within this framework.\smallskip 

  In Section 4, it is shown how two commuting shifts of relative depth two 
and hence multiplicative unitaries arise naturally in terms of 
ambidextrous Hilbert spaces in a tensor $W^*$--category.\smallskip 

  Section 5 is devoted to studying Hilbert spaces in and endomorphism of 
the Cuntz algebra which are algebraic with respect to the natural 
grading. We show, for example, how the problem of determining 
intertwiners between algebraic endomorphisms can be reduced to a 
purely algebraic problem. These results are used in the final section 
to arrive at the duality results already announced earlier in the 
introduction. The paper concludes with an appendix on braided symmetry.
\smallskip

In this paper we prefer to work with strictly associative tensor products and a 
simple way of achieving this is to use as the underlying Hilbert spaces the 
Hilbert spaces in some fixed von Neumann algebra since these are objects 
in a strict tensor $W^*$--category. We will be concerned here with the 
representation categories of multiplicative unitaries and recall the 
basic definitions from \cite{BS}. If $K$ is such a Hilbert space then a unitary 
$V$ on the tensor square $K^2$ is said to be multiplicative if 
$$V_{12}V_{13}V_{23}=V_{23}V_{12},$$
where we use the usual convention regarding indices and tensor products.
A representation of $V$ on a Hilbert space $H$ is a unitary $W\in(HK,HK)$ 
such that 
$$W_{12}W_{13}V_{23}=V_{23}W_{12},\quad \text{on}\quad HK^2.$$
If $W$ and $W'$ are representations of $V$ on $H$ and $H'$ respectively, we say 
that $T\in (H,H')$ intertwines $W$ and $W'$ and write $T\in(W,W')$ if 
$T\times 1_KW=W'T\times 1_K$. We define the tensor product of $W$ and $W'$ 
to be the representation $W\times W'$ on $HH'$ given by
$W\times W':=W_{13}W'_{23}$. The usual tensor product of intertwiners is 
again an intertwiner and in this way we get a strict tensor $W^*$--category 
${\cal R}(V)$ of representations of $V$. In fact this assertion does not 
depend on $V$ being multiplicative. When it is then $V$ itself is a 
representation of $V$ called the regular representation.\smallskip 

   A corepresentation of $V$ on $H$ is a unitary $W\in(KH,KH)$ such that 
$$V_{12}W_{13}W_{23}=W_{23}V_{12}\quad \text{on}\quad K^2H.$$ 
If $W$ and $W'$ are corepresentations on $H$ and $H'$ respectively, we say 
that $T\in(H,H')$ intertwines $W$ and $W'$ and write $T\in(W,W')$ if 
$1_K\times TW=W'1_K\times T$. The tensor product $W\times W'$ of 
corepresentations is defined by $W\times W':=W_{12}W'_{13}$. Just 
as in the case of representations we get a strict tensor $W^*$--category 
now denoted by ${\cal C}(V)$. If $\vartheta=\vartheta_{K,K}$ denotes the 
flip on $K^2$ then $\vartheta V^*\vartheta$ is again a multiplicative unitary 
and the mapping $W\mapsto \tilde W:=\vartheta_{H,K}W^*\vartheta_{K,H}$ defines  
a 1--1 correspondence between representations of $V$ and corepresentations 
of $\vartheta V^*\vartheta$. However, it does not define an  isomorphism 
of tensor $W^*$--categories since $W\times W'\mapsto \tilde W'_{13}\tilde W_{12}$ 
and so leads to an alternative definition of the tensor product of 
corepresentations. In fact the two expressions for the tensor product 
will be equal if and only if $\vartheta_{W,W'}\in(W\times W',W'\times W)$, 
and this corresponds to the case of a group cf.~\cite{CMP}, Prop.~2.5. 
Thus exchanging the definitions of tensor product corresponds to 
exchanging representations of a multiplicative unitary and 
corepresentations of the dual multiplicative unitary.

\section{Preliminaries on Categories of Hilbert Spaces}

  We begin our considerations with a simple but useful lemma on 
natural transformations in the context of $W^*$--categories. Elementary
results and definitions on $W^*$--categories can be found in \cite{GLR} 
and Lemma 2.1 below is just a slight generalization of Corollary 7.4 in \cite{GLR}.
The notion of direct sum will be used in the Hilbert space sense 
rather than in the purely algebraic sense. Thus $A$ is a direct sum of 
objects $B_i$, $i\in I$, if there are isometries $W_i\in(B_i,A)$ 
such that $\sum_iW_i\circ W_i^*=1_A$, where the convergence is in, say, 
the $s$--topology. In particular if the $B_i=B$ for all $i$, 
the condition amounts to saying that there is a Hilbert space of 
support $1_A$ in $(B,A)$. The isometries $W_i$ form an orthonormal 
basis of such a Hilbert space. An object $B$ has {\it central support one} 
or is a {\it generator} if given any object $A$ there are partial 
isometries $W_{i,A}\in(B,A)$ such that $\sum_iW_{i,A}\circ W_{i,A}^*=1_A$, 
see Proposition 7.3 of \cite{GLR}.\smallskip 

\noindent
{\bf 2.1 Lemma} {\sl Let ${\cal T}$ and ${\cal K}$ be $W^*$--categories
and 
$E$ and $F$ be normal $^*$--functors from ${\cal K}$ to ${\cal T}$. 
Suppose ${\cal K}$ has a object $B$ of central support one. 
Then a natural transformation $t$ from $E$ to $F$ has form 
$$t_A=\sum_iF(W_{i,A})\circ T\circ E(W_{i,A})^*,$$ 
where the sum is taken over partial isometries $W_{i,A}$ of $(B,A)$ 
with $\sum_iW_{i,A}\circ W_{i,A}^*=1_A$ and $T$ is an arbitrary element 
of $(E(B),F(B))$ satisfying the intertwining relation
$TE(S)=F(S)T$, $S\in(B,B)$. $t$ is automatically bounded and $t_B=T$.}
\smallskip 

\noindent 
{\bf Proof.} The sum defining $t_A$ converges in the $s$--topology, say, 
and $\|t_A\|\leq \|T\|$. (Consider a finite sum and use the $C^*$--property 
of the norm.) Noting that $t_B=T$, we conclude that $\|t\|=\|T\|$. 
Pick $Y\in(A,C)$, then 
$$t_C\circ E(Y)=\sum_{i,j}F(W_{i,C})\circ T\circ E(W_{i,C}^*
\circ Y\circ W_{j,A}\circ W_{j,A}^*).$$ 
But $W_{i,B}^*\circ Y\circ W_{j,A}\in(B,B)$ so using the intertwining 
property of $T$, we deduce that $t_CE(Y)=F(Y)t_A$ and we have a bounded 
natural transformation.
Conversely, suppose $t\in(E,F)$ and $W\in(B,A)$ then 
$t_A\circ E(W)=F(W)\circ t_B$ implying that $t$ 
is obtained by the above construction with $T=t_B$.\smallskip 

  Note that if $E(B)=F(B)$ and $B$ is irreducible, then we even have a canonical 
natural transformation $t\in(E,F)$, satisfying $t_{B}=1_{E(B)}$. 
We use the notation $(F,T,E)$ for the natural transformation constructed 
as above from $T\in(E(B),F(B))$ satisfying the intertwining relation. The 
usual operations on natural transformations have simple expressions in this 
notation. Thus composition of natural transformations corresponds to 
composing these symbols in the obvious way: 
$$(F,T,E)\circ(E,S,D)=(F,T\circ S,D).$$ 
In fact, the map $t\mapsto t_B$ is a full and faithful $^*$--functor from 
the category of normal $^*$--functors from ${\cal K}$ to ${\cal T}$ to 
the category of normal representations of $(B,B)$. 
If $G$ is a normal $^*$--functor from ${\cal T}$, then acting on $(F,T,E)$ 
on the left with $G$ gives $(GF,G(T),GE)$. If $D$ is a normal 
$^*$--endofunctor of ${\cal K}$ then acting on $(F,T,E)$ on the 
right by $D$ gives $(FD,S,ED)$, where $S=(F,T,E)_{D(B)}$.\smallskip 

\noindent
{\bf 2.2 Lemma} {\sl With the above notation, suppose that $GE=ED$ and $GF=FD$. 
Then if $G(T)=S$, $G\times(F,T,E)=(F,T,E)\times D$. If $G\times(F,T,E)$ 
is invertible, then
$$(F,T,E)\times D=G\times(F,T,E)\circ(ED,G(T)^{-1}\circ S,ED).$$} 

\noindent
{\bf Proof.} As the natural transformations are uniquely determined by 
their values in $B$, the result follows by evaluating in $B$.
\smallskip

  By a category of Hilbert spaces, we mean a $W^*$--category whose objects 
are Hilbert spaces and whose mappings are all bounded linear mappings 
between these Hilbert spaces. Any object in such a category has central 
support one. Any $W^*$--category with an irreducible 
object of central support one is a category of Hilbert spaces in a 
natural way.\smallskip 

  We now consider a category ${\cal K}$ of Hilbert spaces whose objects 
are labelled by ${\Bbb N}_0$ with $(0,0)={\Bbb C}$ and equipped with a
normal 
$^*$--functor $F$ from ${\cal K}$ to ${\cal K}$ such that the action of
$F$ 
on objects is given by $F(n)=n+1$, $n\in{\Bbb N}_0$.\smallskip 

   Since $(0,0)={\Bbb C}$, $(0,n)$ is a Hilbert space of support $1$ for
each 
$n$ and we let $\psi_{i,n}$, $i\in I_n$, be an orthonormal basis of $(0,n)$ 
and set 
$$\hat F(X):=\sum_iF^n(\psi_{i,1})\circ X\circ F^m(\psi_{i,1}^*),\quad X\in(m,n).$$  
Then, with the notation of Lemma 2.1, $\hat F(X)=(F^n,X,F^m)_1$ and $\hat F$ is 
another normal $^*$-functor with $\hat F(n)=n+1$, $n\in{\Bbb 
N}_0$.\smallskip 

 It is obvious that $F$ and $\hat F$ commute. If we iterate 
$\hat F$, we find that 
$$\hat F^k(X)=\sum_iF^n(\psi_{i,k})\circ X\circ
F^m(\psi_{i,k}^*)=(F^n,X,F^m)_k.$$ 

   Since $\hat F$ has the properties assumed for $F$, we can form $\hat{\hat F}$. 
Calculating, we find that 
$$\hat{\hat F}(X)=\sum_{i,j,k}F(\psi_{j,m})\circ\psi_{i,1}\circ\psi_{j,n}^*
\circ X\circ\psi_{k,m}\circ\psi_{i,1}^*\circ F(\psi_{k,m})^*$$ 
$$=\sum_{j,k}F(\psi_{j,m})\circ F(\psi_{j,n}^*\circ X\circ\psi_{k,m})
\circ F(\psi_{k,m})^*=F(X).$$ 
Thus the operation $\hat{}$ is involutive.\smallskip 

  From Lemma 2.1, we know that the natural transformations between powers of 
the functor $F$ are automatically bounded, and furthermore, that a natural 
transformation $t\in(F^r,F^s)$ has the form $t_n=\hat F^n(T)$ with 
$T\in(r,s)$. Finally, we know that $T\mapsto t$ is an isomorphism of 
$W^*$--categories between ${\cal K}$ and the category of natural
transformations 
between the powers of $F$. The latter category is, however, a tensor 
$W^*$--category, so we may use the isomorphism  to equip ${\cal K}$ with 
a tensor product making it into a tensor $W^*$--category. 
We compute this tensor product. If $y\in(F^r,F^s)$ and $y'\in(F^{r'},F^{s'})$, 
then $y\times y'=y\times F^{s'}\circ F^r\times y'$. In other words, 
$(y\times y')_n=y_{n+s'}\circ F^r(y'_n)=\hat F^{n+s'}(Y)\circ F^r\hat F^n(Y')$. 
Setting $n=0$, we see that the tensor product in ${\cal K}$ is given by 
$$Y\times Y':=\hat F^{s'}(Y)\circ F^r(Y').$$ 
Thus $\hat F$ is the functor of tensoring on the right by the object $1$ 
and $F$ the functor of tensoring on the left by the same object.\smallskip 

  The above result can also be seen in a different way: there is 
a functor ${\cal F}$ from ${\cal K}$ into the category of endofunctors of 
${\cal K}$ defined on objects by ${\cal F}(n)=\hat{F}^n$ and on arrows by 
${\cal F}(X)_r:=F^r(X)$, $X\in (m,n)$. It combines the operations. In fact 
if $T$ is any arrow of ${\cal K}$, 
$${\cal F}\times 1_{{\cal F}(1)}={\cal F}F(T),\quad 
1_{{\cal F}(1)}\times{\cal F}(T)={\cal F}\hat{F}(T).$$ 
In view of this result, we refer to a normal $^*$--endofunctor on ${\cal
K}$ 
with $F(n)=n+1$ as being a {\it shift} on ${\cal K}$. 

\smallskip 
  
  As is well known, any two definitions of the tensor product on a category 
of Hilbert spaces are equivalent, i.e.\ the identity functor extends to 
a relaxed tensor functor. To be specific in the case at hand, if we 
have two shifts $F$ and $G$ on ${\cal K}$ and we define 
$V_{m,n}:= (\hat G^n,1_n,\hat F^n)_m$ then 
$$V_{m',n'}\circ\hat F^{n'}(Y)\circ F^m(Z)=\hat G^{n'}(Y)\circ G^m(Z)\circ V_{m,n},$$ 
where $Y\in(m,m')$ and $Z\in(n,n')$. Furthermore, 
$$\hat G^p(V_{m,n})\circ V_{m+n,p}=G^m(V_{n,p})\circ V_{m,n+p}.$$

   If we now define inductively $K:=(0,1)$ and $K^n:=F(K^{n-1})\circ K$, 
where the norm closed linear span is understood, then the above result 
shows that $K^n=(0,n)$.\smallskip 

  We now make some remarks on the automorphisms of ${\cal K}$. 
Such an automorphism will be a normal $^*$--functor $\Gamma$ and we 
suppose, as part of the definition that $\Gamma$ leaves the objects fixed. 
Since an automorphism of a Hilbert space is given by a unitary operator, 
any automorphism of ${\cal K}$ will be inner, meaning that there is a 
unitary natural transformation $u$ from the identity functor to $\Gamma$. 
Lemma 2.1 tells us that $u$ is determined by $u_0\in{\Bbb C}$. Since we
are 
interested in $\Gamma$ rather than $u$, we fix the free phase by requiring 
that $u_0=1$. An inner automorphism is then determined by a sequence of 
unitaries $u_n\in(n,n)$ with $u_0=1$. We now look for the inner 
automorphisms which commute with $F$. They must therefore commute with 
$\hat F$ and preserve the tensor product structure determined by $F$. 
Applying Lemma 2.2, we derive the condition $u_{n+1}=F(u_n)\hat F^n(u_1)$. 
Solving the recurrence relation gives 
 $u_{n+1}=F^n(u_1)F^{n-1}\hat F(u_1)\cdots \hat F^n(u_1)$.
In terms of the tensor product structure determined by $F$, this means 
that the $u_n$ are just tensor powers of $u_1$. If $u_1$ is a phase, 
we get automorphisms which commute with every shift and are the 
analogues of the grading automorphisms of the Cuntz algebra. Of course, 
if $\Gamma$ does not commute with $F$ it maps the tensor product 
structure determined by $F$ onto that determined by $\Gamma F\Gamma^{-1}$.\smallskip

  Our next goal is to characterize all normal $^*$--functors on ${\cal K}$ 
that commute with the given functor $F$. Note first that the action of 
such a functor $G$ on objects must be of the form $G(n)=r+n=F^r(n)$, 
$n\in{\Bbb N}_0$, 
for some $r\in{\Bbb N}_0$. Thus $(G,1_r,F^r)$ will be a natural unitary 
transformation from $F^r$ to $G$ and since $F$ commutes with these two 
functors, we may apply the second identity of Lemma 2.2 to deduce that 
$$R_{n+1}=F(R_n)\circ\hat F^n(R_1),$$ 
where we have written $R$ for $(F^r,1_r,G)$. 

Conversely, given a unitary $R_1\in(r+1,r+1)$, 
take $R_0$ to be the unit on $r$ and define $R_n$, $n>1$, inductively using 
the above formula. Finally, define 
$$G(X)=R_n\circ F^r(X)\circ R_m^*.$$ 
Then $G$ is obviously a normal $^*$--functor with $G(n)=n+r$, $n\in{\Bbb 
N}_0$. 
Furthermore, if $X\in(m,n)$,  
$$GF(X)=R_{n+1}\circ F^{r+1}(X)\circ R_{m+1}^*
=F(R_n)\circ\hat F^n(R_1)\circ F^{r+1}(X)\circ R_{m+1}^*$$ 
$$=F(R_n)\circ F^{r+1}(X)\circ F(R_m)^*=FG(X).$$
Thus $F$ and $G$ commute and we have proved the following result.\smallskip  

\noindent 
{\bf 2.3 Proposition} {\sl  
Normal $^*$--functors 
$G$ commuting with $F$ on ${\cal K}$ 
are of the form
$$G(\psi)=R\circ F^r(\psi), \quad \psi\in(0,1),\  R\in(r+1, r+1).$$ 
$r$ is called the rank of $G$.}\smallskip 

  Note that since we have already computed all natural transformations 
between the tensor powers of $F$, we have implicitly computed all natural transformations 
betwen two functors $G$ and $G'$ commuting with $F$. Note, too, 
that $G$ may be obtained from $F^r$ by acting on the left with an inner 
automorphism.\smallskip 

It is easy to compute what composition of functors means for the
corresponding unitary operators. If we use the notation $G_R$ to denote
the 
functor corresponding to the unitary $R\in(r+1,r+1)$, then $G_RG_S=G_T$,
where 
$T=G_R(S)\circ F^s(R)$. Note that Proposition 2.2 is closely related to Cuntz's 
result\cite{Cu1} characterizing the endomorphisms of the Cuntz algebra in terms of 
unitary operators.\smallskip

  We now consider another category ${\cal H}$ of Hilbert spaces whose 
objects will be denoted $H_n$ with $n\in{\Bbb N}_0$. Suppose we have a
normal 
$^*$--functor $H$ from ${\cal K}$ to ${\cal H}$ whose action on objects 
takes $n$ to $H_n$. Then we may define a normal $^*$-endofunctor $A$ on 
${\cal H}$ by setting 
$$A(X):=\sum_iH\hat F^n(\psi_{i,1})\circ X\circ H\hat F^m(\psi_{i,1})^*
=(H\hat F^n,X,H\hat F^m)_1,\quad X\in(H_m,H_n),$$ 
where, as before, the sum runs over an orthonormal basis. We see at 
once that, with this definition, $AH=HF$. $H$ is to be thought of as 
tensoring on the left by $H_0$ and $A$ as tensoring on the right by 
the object $1$ of ${\cal K}$. An easy calculation now shows that 
$$A^s(X)\circ H\hat F^m(Y)=H\hat F^n(Y)\circ A^r(X),\quad X\in(H_m,H_n),\,Y\in(r,s).$$ 
These equations show that we may define a functor ${\cal H}$ from 
${\cal K}$ to End${\cal H}$ by setting ${\cal H}(r)=A^r$ and 
${\cal H}(Y)_n=H\hat F^n(Y)$ for $Y\in(r,s)$.\smallskip 

\noindent
{\bf 2.4 Lemma} {\sl There is a $1-1$ correspondence between normal 
$^*$--functors $H$ from ${\cal K}$ to ${\cal H}$ with $H(n)=H_n$ and 
normal $^*$--functors ${\cal H}$ from ${\cal K}$ to End${\cal H}$ such 
that ${\cal H}(Y)\times\hat F={\cal H}\hat  F(Y)$ given 
by $H(Y):={\cal H}(Y)_0$.}\smallskip 

\noindent 
{\bf Proof.} We have seen above how to construct the functor ${\cal H}$
from 
$H$. 
$$({\cal H}(Y)\times F)_r={\cal H}(Y)_{r+1}=HF^{n+1}(Y)=({\cal
H}F(Y))_r.$$ 
Conversely, given ${\cal H}$, we obtain $H$ by evaluating in $0$: 
$H(Y):={\cal H}(Y)_0$. Now $HF^n(Y)={\cal H}(F^n(Y))_0={\cal
H}(Y)_n$. Thus 
${\cal H}$ is the functor associated with $H$.\smallskip

  We have seen how ${\cal K}$ with the functor $F$ is isomorphic to the 
category of tensor powers of $(0,1)$ such that $F$ becomes the functor of 
tensoring on the left by the object $1$. There is a similar result for 
${\cal H}$ with the functor $H$.\smallskip 

\noindent
{\bf 2.5 Proposition} {\sl Let ${\cal H}$ be a category of Hilbert spaces 
and $H:{\cal K}\to{\cal H}$ a normal $^*$--functor with $H(r)=H_r$. Then 
there is a unique isomorphism $\Phi$ of $(H_0,H_0)\otimes{\cal K}$ 
into ${\cal H}$ such that $\Phi(T\otimes Y)=A^s(T)H(Y)$,
$Y\in(r,s)$.}\smallskip 

\noindent
{\bf Proof.} Noting that $A^s(T)H(Y)=H(Y)A^r(T)$, we see that $\Phi$ extends 
uniquely to the algebraic tensor product. Now given $X\in(H_r,H_s)$, 
write $$X=\sum_{i,j}H(\psi_{i,s})\circ X_{ij}\circ H(\psi_{j,r})^*,$$ 
where $X_{ij}:=H(\psi_{i,s})^*\circ X\circ H(\psi_{j,r})\in(H_0,H_0)$. 
Hence $X=\sum_{i,j}A^s(X_{ij})\circ H(\psi_{i,s}\circ\psi_{j,r}^*)$  
showing that $\Phi$ extends by continuity to a full functor. But any 
normal $^*$--functor from ${\cal K}$ is faithful hence $\Phi$ is  
an isomorphism.\smallskip 

  In view of the above results, given a shift $F$ on ${\cal K}$, we 
may say that a normal $^*$--functor $H$ from ${\cal K}$ to ${\cal H}$ 
with $H(r)=H_r$ determines an {\it action} of $({\cal K},F)$ 
on ${\cal H}$ via Lemma 2.4 making ${\cal H}$ into a (right) 
$({\cal K},F)$--module.
\smallskip 

  We clearly should be able to define the tensor product of two 
${\cal K}$--modules ${\cal H}$ and ${\cal H}'$ and it should just involve 
replacing $(H_0,H_0)$ and $(H'_0,H'_0)$ by $(H_0,H_0)\otimes (H'_0,H'_0)$. 
Let $H$ and $H'$ be actions of ${\cal K}$ on ${\cal H}$ and ${\cal H}'$, 
respectively. Then the tensor product of the two actions is an action 
$H\otimes H'$ on a category of Hilbert spaces ${\cal H}\otimes{\cal H}'$ 
together with normal $^*$--functors $D:{\cal H}\to{\cal H}\otimes{\cal
H}'$ 
and $D':{\cal H}'\to{\cal H}\otimes{\cal H}'$ such that $H\otimes
H'=DH=D'H'$ 
and $DA=A\otimes A'D$ and $D'A'=A\otimes A'D'$, where $A\otimes A'$ is the 
endofunctor on ${\cal H}\otimes{\cal H}'$ associated with $H\otimes H'$. 
In restriction to $(H_0,H_0)$ and $(H'_0,H'_0)$, we require that $D$ and $D'$ 
should define $((H\otimes H')_0,(H\otimes H')_0)$ as a tensor product of 
the von Neumann algebras $(H_0,H_0)$ and $(H'_0,H'_0)$. It should be noted 
that $D$ and $D'$ are uniquely determined by their values on $(H_0,H_0)$ 
and $(H'_0,H'_0)$, respectively.  
In fact, we have seen that with a suitable definition of tensor product the 
action $H$ becomes $1_{(H_0,H_0)}\otimes:{\cal
K}\to(H_0,H_0)\otimes{\cal K}$ 
with a similar expression for $H'$. Taking 
${\cal H}\otimes{\cal H'}=(H_0,H_0)\otimes(H'_0,H'_0)\otimes{\cal K}$ and 
$H\otimes H'(Y):=1_{(H_0,H_0)\otimes(H'_0,H'_0)}\otimes Y$ and defining $D$ 
to be the normal $^*$--functor such that 
$D(T\otimes Y):=T\otimes 1_{(H'_0,H'_0)}\otimes Y$, for $T\in(H_0,H_0)$ and 
$Y\in(r,s)$, with a similar expression for $D'$ we do get a tensor product 
of $H$ and $H'$.\smallskip  

  In fact, every tensor product of actions is of this form, since, 
writing $H''$ for $H\otimes H'$, we may use $D$ and $D'$ to realize 
$(H''_0,H''_0)$ as a tensor product of $(H_0,H_0)$ and $(H'_0,H'_0)$. 
Then since $H''$ is an action, we have an isomorphism $\Phi''$ from 
$H''\otimes{\cal K}$ to ${\cal H}\otimes{\cal H}'$ and $\Phi$ and $\Phi'$ 
from $H\otimes{\cal K}$ and $H'\otimes{\cal K}$ to ${\cal H}$ and 
${\cal H}'$, respectively. If we then define 
$d(T\otimes Y):=T\otimes 1_{(H'_0,H'_0)}\otimes Y$ and 
$d'(T'\otimes Y):=1_{(H_0,H_0)}\otimes Y$ using the tensor product 
structure on $(H''_0,H''_0)$ coming from $D$ and $D'$, we have 
$\Phi''\, d=D\Phi$ and $\Phi''\, d'=D'\Phi$. Thus the action $H''$ is 
ismorphic to an explicit tensor product of actions.\smallskip 

We now look for generalizations of our result characterizing 
categories ${\cal K}$ of Hilbert spaces with $(0,0)={\Bbb C}$ 
which admit a shift by dropping the condition that $0$ should 
be irreducible.\smallskip 

\noindent
{\bf 2.6 Lemma} {\sl Let $H$ and $K$ be Hilbert spaces and $F$ a 
unital normal homomorphism from $(H,H)$ to $(K,K)$. Then we can pick 
matrix units $E_{ab}$ in $(H,H)$ and $E_{ar,bs}$ in $(K,K)$ 
such that $F(E_{bc})=\sum_rE_{br,cr}$.}\smallskip 

\noindent
{\bf Proof.} Pick matrix units $E_{ab}$ in $(H,H)$ then for a fixed 
$a$, $F(E_{aa})$ is a non-zero projection and we may pick matrix 
units $E_{ar,as}$ for the Hilbert space $F(  E_{aa})K$. Thus, in particular, 
$F(E_{aa})=\sum_rE_{ar,ar}$. Now define
$$E_{br,cs}:=F(E_{ba})E_{ar,as}F(E_{ac}).$$ 
A routine computation shows that 
$$E_{br,cs}E_{dt,ev}=\delta_{cd}\delta_{st}E_{br,ev}$$ 
$$\sum_rE_{br,cr}=F(E_{bc})=I.$$ 
Thus we have defined a set of matrix units with the required properties.\smallskip 

\noindent
{\bf 2.7 Lemma} {\sl Let $H$ and $K$ be Hilbert spaces and $F$ a unital 
normal homomorphism from $(H,H)$ to $(K,K)$ then there is a Hilbert space $L$ 
of support one in $(H,K)$ such that 
$$\psi T=F(T)\psi,\quad \psi\in L,\,\,T\in(H,H).$$}

\noindent
{\bf Proof.} Pick orthonormal bases $e_a$ and $e_{ar}$ in $H$ and $K$ 
in such a way that the corresponding matrix units are as in Lemma 2.6. 
Define $\psi_r\in(H,K)$ by $\psi_re_a:=e_{ar}$ and a computation 
shows that $\psi_r$ is a basis of a Hilbert space $L$ of support 
one in $(H,K)$. Now $\psi_rE_{bc}e_d=\delta_{cd}e_{br}$ and 
$F(E_{bc})\psi_r e_d=\sum_sE_{bs,cs}e_{dr}=\delta_{cd}e_{br}$, 
completing the proof.\smallskip

  Now let ${\cal H}$ be a category of Hilbert spaces 
with objects labelled by ${\Bbb N}_0$ and $F$ a 
normal $^*$-functor from ${\cal H}$ to ${\cal H}$ acting on objects 
as $F(n)=n+1$. We have analyzed the situation 
when $0$ is irreducible and we begin a similar analysis in the general 
case.\smallskip 

   By Lemma 2.7, there is a Hilbert space $L$ of support one 
in $(0,1)$ inducing $F$ on $(0,0)$. We now define as before 
$$\hat F(X)=\sum_iF^n(\psi_i)\circ X\circ F^m(\psi^*_i).$$ 
$\hat F$ will be a normal $^*$--functor with $\hat{F}(n)=n+1$ and
commuting 
with $F$ and $\hat{\hat F}= F$. The
$$F^{n-1}(\psi_{i_1})\circ\cdots\circ F(\psi_{n-1})\circ \psi_{i_n}$$ 
form an orthonormal basis $\psi_{i,n}$ of a Hilbert space of support 
one in $(0,n)$. We find 
$$\hat F^k(X)=\sum_iF^n(\psi_{i,k})\circ X\circ F^m(\psi^*_{i,k}).$$ 
We can interchange the role of $F$ and $\hat F$. But care is required 
as we need another Hilbert space of support one in $(0,n)$. It follows 
from Lemma 2.7 that $\hat F(T)=F(T)$ for $T\in(0,0)$ so that $L$ 
remains the correct Hilbert space in $(0,1)$. In $(0,n)$ we need a 
Hilbert space with orthonormal basis 
$$\hat F^{n-1}(\psi_{i_1})\circ\cdots\circ\hat F(\psi_{n-1})\circ \psi_{i_n}$$ 
denoted $\hat\psi_{i,n}$. Then 
$$F^k(X)=\sum_i\hat F^n(\hat\psi_{i,k})\circ X\circ\hat F^m(\hat\psi^*_{i,k}).$$ 

  We now define a $W^*$--subcategory ${\cal K}$ of ${\cal H}$ whose arrows 
consist of those $S\in(m,n)$ such that $S\circ F^m(T)=F^n(T)\circ S$. 
This subcategory is invariant under $F$ and the object $0$ is now irreducible. 
By Lemma 2.7, $L$ is a Hilbert space of support one in this category 
so that this is actually a category of Hilbert spaces. So our previous 
result is applicable and ${\cal K}$ is isomorphic to the category of 
tensor powers of a Hilbert space $K$ in such a way that $F$ is the 
functor of tensoring on the right by $1_K$.\smallskip 
  
  But the inclusion functor of ${\cal K}$ in ${\cal H}$ is now a normal 
$^*$--functor which is the identity on objects. The associated endofunctor 
$A$ is just $F$. Our previous analysis yields the following result.\smallskip 

\noindent
{\bf 2.8 Proposition} {\sl Let ${\cal H}$ be a category of Hilbert spaces 
with objects labelled by ${\Bbb N}_0$ and $F$ a normal $^*$-functor from 
${\cal H}$ to ${\cal H}$ acting on objects as $F(n)=n+1$. 
Then ${\cal H}$ is isomorphic to the category of Hilbert spaces whose 
objects are of the form $H\otimes K^n$, $n\in{\Bbb N}_0$, for Hilbert 
spaces of $H$ and $K$ in such a way that $F$ becomes the functor of tensoring 
on the right by $1_K$.}\smallskip

\section{Shifts and Multiplicative Unitaries}  

We used above the characterization of endomorphisms of the 
Cuntz algebra in terms of unitaries given by Cuntz who also
showed\cite{Cu2} 
that if $K$ is a finite-dimensional Hilbert 
space then a unitary $V\in(K^2,K^2)$ is a multiplicative unitary if the endomorphism $\lambda_R$ 
of the Cuntz algebra ${\cal O}_K$ satisfies 
$$\lambda_R\lambda_R=\rho\lambda_R,$$ 
where $R=V\theta$, $\theta$ being the flip on $K^2$, and $\rho$ the endomorphism 
generated by the defining Hilbert space $K$. This result remains valid 
in the extended Cuntz algebra if $K$ is infinite dimensional\cite{R}. The 
multiplicativity of $V$ is alternatively expressed by the identity 
$\lambda_{V^*}\lambda_R=\rho\lambda_{V^*}$. The corresponding identities 
for normal $^*$-functors on ${\cal K}$ commuting with $F$ imply that 
a unitary $V\in(2,2)$ is multiplicative for the tensor structure 
induced by $\hat F$, i.e.\ where $F$ corresponds to tensoring on the right 
and $\hat F$ to tensoring on the left by the object $1$. The role of the 
endomorphism $\rho$ of the Cuntz algebra is played by $\hat F=G_\theta$. 
In other words, the following result holds.\smallskip 

\noindent
{\bf 3.1 Proposition} {\sl Let $F$ and $G$ be two commuting shifts 
on ${\cal K}$. Let $R\in(2,2)$ be the unitary such that 
$G(\psi)=R\circ F(\psi)$, $\psi\in(0,1)$. Set $V:=R\circ\theta$, 
where $\theta$ is the flip on $(2,2)$ derived from the tensor structure 
induced by $\hat F$, then the following conditions are equivalent. 
\begin{description}
\item{a)} $V$ is a multiplicative unitary.
\item{b)} $GG=\hat FG$.
\item{b')} $FF=\hat GF$.
\item{c)} $n\mapsto F^{n+1}(\psi)$ is a natural transformation from 
$G$ to $GG$, $\psi\in(0,1)$.
\item{c')} $n\mapsto G^{n+1}(\psi)$ is a natural transformation from 
$F$ to $FF$, $\psi\in(0,1)$.
\item{d)} $GG(\psi)\circ F(\psi')=FF(\psi')\circ 
G(\psi)$, $\psi,\psi'\in(0,1)$.
\item{e)} $G_{V^*}G=\hat FG_{V^*}$. 
\item{f)} $n\mapsto F^{n+1}(\psi)$ is a natural transformation from 
$G_{V^*}$ to $G_{V^*}G$, $\psi\in(0,1)$.
\item{g)} $G_{V^*}G(\psi)\circ F(\psi')=FF(\psi')\circ G_{V^*}(\psi)$, $\psi,\psi'\in(0,1)$.
\end{description}} 

\noindent
{\bf Proof.} The equivalence of a), b) and e) is a simple computation 
whose origins were explained above. Suppose a) is valid and pick $X\in(m,n)$ 
then 
$$\sum_iF^{n+1}(\psi_{i,1})\circ G(X)\circ F^{m+1}(\psi_{i,1})^*=\hat FG(X)=GG(X).$$
Hence $F^{n+1}(\psi)\circ G(X)=G(X)\circ F^{m+1}(\psi)$ for $\psi\in(0,1)$, 
giving c). d) follows as a special case. But if d) holds, then 
$$\hat FG(\psi)=\sum_iF^2(\psi_{i,1})G(\psi)F(\psi_{i,1})^*=GG(\psi).$$ 
But the set of arrows $X$ in ${\cal K}$ such that $\hat FG(X)=GG(X)$, is a 
$W^*$--subcategory of ${\cal K}$ which is invariant under the action of
$F$, 
seeing that $F$ commutes with $\hat F$ and $G$. Thus d) implies b). In the 
same way, we can prove that e), f) and g) are equivalent. Finally, the 
symmetry between $F$ and $G$, visible in d), shows that b') and c') are 
also equivalent to the remaining conditions.\smallskip

  We have seen that the situation in Proposition~3.1 is symmetric in the two 
commuting shifts $F$ and $G$. Interchanging $F$ and $G$ obviously corresponds 
to replacing $R$ by $R^{-1}$ and hence  $V:=R\circ\theta$ by 
$\hat V=R^{-1}\circ\theta$, the dual multiplicative unitary. The multiplicative 
unitaries on $(0,2)$ for the tensor structure determined by $F$ are 
clearly in $1-1$ correspondence with the shifts commuting with $F$ 
and satisfying the equivalent conditions of Proposition~3.1. Equivalent 
multiplicative unitaries correspond to shifts conjugate under 
(inner) automorphisms of ${\cal K}$ commuting with $F$.\smallskip 

 We can also characterize the intertwining operators between the tensor 
powers of $V$, regarded as an object in the category of representations of 
$V$, in terms of the commuting shifts.\smallskip 

\noindent
{\bf 3.2 Lemma} {\sl Let $F$ and $G$ be commuting shifts on ${\cal K}$
with 
$GG=\hat FG$ and let $V$ be the associated multiplicative unitary. Set 
$G':=G_{\vartheta V}$ then 
$$(V^r,V^s)=\{ Y\in(r,s):G'(Y)=\hat F(Y)\}.$$} 

\noindent
{\bf Proof.} Interpreting $F$ as tensoring on the right by $1$, the 
condition $G'(Y)=\hat F(Y)$ reads 
$$(\vartheta V)_s\circ Y\times 1\circ (\vartheta V)_r^*
=\vartheta_s\circ Y\times 1\circ \vartheta_r^*.$$ 
Hence, it suffices to show 
that $\vartheta_s^*\circ (\vartheta V)_s=V_{1s+1}V_{2s+1}\cdots V_{ss+1}$, 
where we have used the index notation on the right hand side. However, 
this may be proved by induction. In fact 
$$\vartheta_s^*\circ(\vartheta V)_s=1_{s-1}\times\vartheta\circ\vartheta_{s-1}^*
\times 1\circ(\vartheta V)_{s-1}\times 1\circ1_{s-1}\times(\vartheta V)$$
$$=\vartheta_{ss+1}V_{1s}V_{2s}\cdots V_{s-1s}\vartheta_{ss+1}V_{ss+1}=
V_{1s+1}V_{2s+1}\cdots V_{ss+1}.$$ 

\noindent
{\bf 3.3 Corollary} {\sl Let $F$ and $G$ be two commuting shifts on ${\cal
K}$ 
such that $GG=\hat FG$ and $V$ the corresponding multiplicative unitary. 
Let $G':=G_{\vartheta V}$ and consider a sequence $t_n\in(r+n,s+n)$ that 
defines simultaneously a natural transformation $t$ in 
$(\hat G^{'r},\hat G^{'s})$ and $(\hat F^r,\hat F^s)$. 
Such natural transformations form a tensor $W^*$--subcategory of the
tensor $W^*$--category of all  natural transformations between the 
powers of $\hat G'$. Evaluating $t$ in $0$ establishes an isomorphism with 
the tensor $W^*$--category of intertwiners between powers of $V$. 
}\smallskip 

\noindent
{\bf Proof.} If $t\in(\hat G^{'r},\hat G^{'s})\cap(\hat F^r,\hat F^s)$, 
then, by Lemma 2.1, $t_n=G^{'n}(t_0)=F^n(t_0)$. Thus $t$ 
is uniquely determined by $t_0$ and, as $G'$ and $F$ commute, the only 
constraint on $t_0$ is that $G'(t_0)=F(t_0)$. Lemma 3.2 shows 
that we have an isomorphism of $W^*$--categories and a computation shows 
that it preserves tensor products.\smallskip

  We now comment on the role of multiplicative unitaries,
seen in this light. As we have seen,   $F$ determines a tensor product structure 
on ${\cal K}$ and $G$ determines another so it is natural to interpret 
$R$ as describing the transition from one tensor product to another. However, 
any unitary $R\in(2,2)$ would serve here. We do not need $V:=R\theta$ to be a 
multiplicative unitary. Although it is perfectly correct to say that 
$F$ determines a tensor product structure on ${\cal K}$, this structure 
really involves two commuting shifts $F$ and $\hat F$ which we interpret 
as tensoring on the right by $1$ and on the left by $1$. Thus  
taking two commuting shifts $F$ and $G$ can be regarded as 
generalizing  the idea of a tensor product. It is less symmetric in 
that $\theta$ has been replaced by $R$ and $\hat F$ by $G$ and does 
not give rise to a bifunctor as a true tensor product unless $G=\hat F$. Looked 
at this way, the condition of $V$ being multiplicative being equivalent to 
b) has a simple interpretation. It is a `depth 2' condition: it is 
not necessary to apply $G$ more than once since in successive applications 
it can always be replaced by $\hat F$. The justification for adopting this 
terminology from the theory of subfactors is Lemma 6.3 of \cite{L}. 
We shall say that the two commuting shifts have relative depth two 
to stress that the notion involves the two shifts 
symmetrically.\smallskip

   Let us look at some examples of actions of ${\cal K}$ on ${\cal H}$. 
A first example is suggested directly by Proposition 2.5. Given a Hilbert 
space $H$, take ${\cal H}:=H\otimes{\cal K}$ and then define 
$H:{\cal K}\to{\cal H}$ by $H(Y):=1_H\otimes Y$. With this definition
we find $A(X)=X\otimes 1_K$. We get other normal $^*$--functors 
from ${\cal K}$ to ${\cal H}$ by picking for each $r$ a unitary 
$W_r\in(H_r,H_r)$ and then defining $E(Y):=W_s\circ 1_H\otimes Y\circ W_r^*$,
$Y\in(r,s)$. 
In fact, we are in the situation of Lemma 2.1 and $W_s=(E,1_{H_0},H)_s$. 
A computation shows that in this case the corresponding endofunctor $B$, say, 
is given by 
$$B(X)=W_{n+1}\circ A(W_n^*\circ X\circ W_m)\circ
W_{m+1}^*,\,\,X\in(H_m,H_n).$$

  In particular, suppose we start from a multiplicative unitary $V$ on $K^2$ 
and a representation $W$ on $H$. We can define $H$ and $E$ as above but 
making the particular choice 
$$W_s:=W_{12}W_{13}\dots W_{1s+1},$$ 
using index notation.\smallskip 

\noindent
{\bf 3.4 Lemma} {\sl Let $W$ be a representation of a multiplicative unitary 
then defining functors $E,H:{\cal K}\to{\cal H}$, as above and letting 
$G$ be the endofunctor on ${\cal K}$ defined by $R:=V\theta$, as before,
we 
have $EG=HG$.}\smallskip 

\noindent 
{\bf Proof.} It suffices to show that $EG(\psi)=HG(\psi)$ for each 
$\psi\in(0,1)$. Now $HG(\psi)=H(R)HF(\psi)$, whereas 
$$EG(\psi)=E(R)\circ EF(\psi)=E(R)\circ W_{12}W_{13}HF(\psi)W_{12}^*
=E(R)\circ W_{12}HF(\psi).$$ 
However $E(R)=W_{12}W_{13}R_{23}W_{13}^*W_{12}^*$ and $H(R)=R_{23}$, 
so the identity in question follows from the definition of a representation 
$W$ of $V$, $W_{12}W_{13}V_{23}=V_{23}W_{12}$.\smallskip 

  We can therefore adopt the following point of view. Regard a 
multiplicative unitary as a category ${\cal K}$ of Hilbert spaces, 
as above, equipped with two commuting shifts 
$F$ and $G$ satisfying $GG=\hat FG$. A representation 
of such a multiplicative unitary is a category of Hilbert spaces 
${\cal H}$, as above, equipped with two normal $^*$--functors $E$ and $H$ 
from ${\cal K}$ to ${\cal H}$ of rank zero such that $EG=HG$. An
intertwining 
operator between two representations is a bounded linear operator 
$T\in(H_0,H'_0)$ such that, setting 
$$A(X):=\sum_iH'\hat F^n(\psi_{i,1})\circ X\circ H\hat F^m(\psi_{i,1})^*,
\quad X\in(H_m,H'_n),$$ 
$r\mapsto A^r(T)$ defines a natural transformation from $E$ to $E'$. In fact, 
if we pick $\psi\in(0,1)$, this gives 
$$A(T)W1_H\otimes \psi=W'1_{H'}\otimes\psi T=W'A(T)1_{H'}\otimes\psi.$$ 
Since $A(T)=T\otimes 1_K$, $T\in(W,W')$. On the other hand, if $T\in(W,W')$, 
then $A^r(T)W_r=W'_rA^r(T)$ and this implies that $r\mapsto A^r(T)$ is a 
natural transformation from $E$ to $E'$. Note, however, that 
$A^r(T)=(H',T,H)_r$, so that $T\in(H_0,H'_0)$ is in $(W,W')$ if 
and only if $r\mapsto A^r(T)$  defines an element of
$(H,H')\cap(E,E')$.\smallskip 

  Superficially, the relation $GG=\hat FG$ looks like a special case 
of $EG=HG$. However, $E$ and $H$ have rank zero, whereas $G$ and 
$\hat F$ have rank one. However, if we delete the object $0$ from 
${\cal K}$ to give a full subcategory ${\cal K}_+$ and let $G_+$ and 
$\hat F_+$ denote the functors obtained by restricting the range to 
${\cal K}_+$, the relation $G_+G=\hat F_+G$ is a special case of 
$EG=HG$ and we are considering $V$ as an object in its category of 
representations.\smallskip 

  The relation $EG=HG$ is symmetrical in $E$ and $H$. Exchanging $E$ 
and $H$ corresponds to replacing $W$ by $W^{-1}$ which will not, in 
general be a representation of the multiplicative unitary $V$. However, 
the tensor product notation involved in the definition of representation 
refers to $H$.\smallskip 
   
   $H$ is to be interpreted here as the trivial representation on 
the same space $H_0$ as $E$. 
Since $(E, E')\cap(H, H')\subset(H, H')$,
the category of representations of $V$ automatically 
comes equipped with a faithful functor into the subcategory of 
trivial representations.\smallskip

 When we have two commuting shifts and two actions $E$ and $H$, as above, 
it seems appropriate to refer to ${\cal H}$ as a 
${\cal K}$--bimodule. We have already considered the tensor product 
of ${\cal K}$--modules in the last section and we now extend these 
considerations to bimodules. Given therefore two action $E'$ and $H'$ 
on ${\cal H}'$, we form ${\cal H}\otimes{\cal H}'$ and $H\otimes H'$. 
Of course, we could equally well have formed $E\otimes E'$ instead, 
but what we really need is the relation between the two actions. 
We therefore use the functors $D$ and $D'$, expressing ${\cal
H}\otimes{\cal H}'$ 
as a tensor product and set 
$$(E\otimes E')(Y)=D(W_s)\circ D'E'(Y)\circ D(W_r^*),\quad Y\in(r,s)$$ 
where $W_r=(E,1_{H_0},H)_r$. It is easy to check that $E\otimes E'$ is 
a normal $^*$--functor and that this definition applied to $H$ and $H'$ 
gives $H\otimes H'$. Expressing $E'$ in terms of $H'$ using 
$W'_r=(E',1_{H'_0},H')_r$, the definition is equivalent to 
$$(E\otimes E',1_{(H\otimes H')_0},H\otimes H')_r=D(W_r)\circ D'(W'_r).$$ 
There is obviously a convention involved here because we have chosen to 
write the primed terms to the right of the unprimed terms. However, this 
convention is consistent with that used for multiplicative unitaries 
in that it corresponds to taking the tensor product of representations 
$W$ and $W'$ as $W_{13}W'_{23}$. It is easily checked that 
$EG=HG$ and $E'G=H'G$ imply $E\otimes E'G=H\otimes H'G$.\smallskip 

\smallskip

\section{Multiplicative Unitaries and Tensor Categories}  

  We now exhibit a mechanism leading from objects in a tensor 
$W^*$--category to multiplicative unitaries. It will involve a 
category of Hilbert spaces ${\cal K}$, equipped with two commuting shifts
$F$ and $G$. Let $\rho$ be an object in a strict tensor 
$W^*$--category and suppose that $K$ is a Hilbert space of support 
one contained in $(\rho,\rho^2)$. We then define inductively 
$K^n:=K^{n-1}\times 1_\rho\circ K$, where the norm closed linear span 
is understood. Then $K^n$ is a Hilbert space of support 1 in 
$(\rho,\rho^{n+1})$. $K^0$ will denote ${\Bbb C}1_\rho$. We now set 
$$(K^m,K^n):=\{X\in(\rho^{m+1},\rho^{n+1}):X\circ K^m\subset K^n\}.$$ 
We see that $K^n=(K^0,K^n)$ and thus we have defined a 
$W^*$--subcategory ${\cal K}$ of Hilbert spaces of the tensor
$W^*$--category 
${\cal T}_\rho$ whose objects are the tensor powers of $\rho$. We claim 
that ${\cal K}$ is invariant under tensoring on the right by $1_\rho$. 
In fact, it suffices to show that $K\times 1_\rho\subset (K,K^2)$, 
i.e.\ that $K\times 1_\rho\circ K\subset K^2$ but this is true by 
construction. If $F$ denotes the restriction of $\times 1_\rho$ to ${\cal
K}$, 
${\cal K}$ has unique structure of tensor $W^*$--category such that $F$ 
becomes the functor of tensoring on the right by $1_K$.\smallskip 

\noindent
{\bf 4.1 Lemma} {\sl Let $K\subset(\rho,\rho^2)$ be a Hilbert space of support 
one in a tensor $W^*$--category and let ${\cal K}$ be the subcategory of
Hilbert 
spaces defined as above, then the following conditions are equivalent.
\begin{description}
\item{a)} $1_\rho\times K\subset (K,K^2)$.
\item{b)} ${\cal K}$ is an invariant subcategory for $1_\rho\times$. 
\item{c)} $K^n=1_\rho\times K^{n-1}\circ K$, $n\in{\Bbb N}$. 
\item{d)} $K^2=1_\rho\times K\circ K$. 
\end{description}} 

\noindent
{\bf Proof.} ${\cal K}$ is the smallest $W^*$--subcategory containing $K$
and 
invariant under $\times 1_\rho$. Furthermore $1_\rho\times$ and $\times 1_\rho$ 
commute, thus a) implies b). Given b), we know from Proposition 2.3 that 
$1_\rho\times\psi=R_n\circ\psi\times 1_\rho$, for $\psi\in K^n$ with $R_n$ 
unitary, giving c). c) implies d), trivially and d) implies a).\smallskip 

  We call a Hilbert space $K$ satisfying the equivalent conditions of 
Lemma~4.1 {\it ambidextrous}.\smallskip 

\noindent
{\bf 4.2 Theorem} {\sl Let $K\subset(\rho,\rho^2)$ be an ambidextrous Hilbert 
space of support one in $(\rho,\rho^2)$. Let $F$ and $G$ denote the restrictions 
of $\times 1_\rho$ and $1_\rho\times$ to ${\cal K}$ then there is a unique 
$V\in(K^2,K^2)$ such that 
$$G(\psi)=V\circ\hat F(\psi),\quad \psi\in K.$$ 
$V$ is a multiplicative unitary.}\smallskip 

\noindent
{\bf Proof.} $V$ is unique and is unitary because it is given by 
$$V=\sum_iG(\psi_{i,1})\hat F(\psi_{i,1})^*$$ 
where the sum is taken over an orthonormal basis. Now since 
$K\in(\rho,\rho^2)$, for $X\in(K^m, K^n),$ 
$$\hat FG(X)=\sum_iF^{n+1}(\psi_{i,1})\circ
 1_\rho\times X\circ F^{m+1}(\psi_{i,1})^*
=1_{\rho^2}\times X=G^2(X).$$
Thus $G^2=\hat FG$ and $V$ is a multiplicative unitary by Proposition 2.3.\smallskip 

  Every multiplicative unitary can be realized in this manner.\smallskip 

\noindent
{\bf 4.3 Proposition} {\sl Let $V$ be a multiplicative unitary $V\in(K^2,K^2)$, 
then $G_{V^*}(K)$ is an ambidextrous Hilbert space $H\subset(K,K^2)$ of 
support one and if $U\in(K,H)$ is the unitary taking $\psi$ to 
$G_{V^*}(\psi)$ the multiplicative unitary defined by 
$H$ is $U\times U\circ V\circ (U\times U)^{-1}$.}\smallskip 

\noindent
{\bf Proof.} $H:=V^*\circ K\times 1_K=G_{V^*}(K)$ is 
obviously a Hilbert space of support one in $(K,K^2)$. To verify that 
$H$ is ambidextrous, it suffices by Proposition 2.4 to verify that 
$$\psi_2^*\times 1_K\circ V\circ 1_K\times\psi_1^*\circ 1_K
\times V\circ V^*\times 1_K\circ \psi_3\times 1_{K^2}\circ V^*\circ\psi_4\times 1_K$$
is a multiple of $1_K$ for any choice of $\psi_i\in K$, $i=1,2,3,4$. 
Using the multiplicativity of $V$, this reduces to 
$$(\psi_2^*\times\psi_1^*\circ V^*\circ\psi_3\times\psi_4)\times 1_K$$ 
and hence is a multiple of $1_K$. It remains to compute the operator 
$S$ defined by 
$$S\circ V^*\times 1_K\circ\psi\times 1_{K^2}=1_K\times V^*\circ 
1_K\times\psi\times 1_K,\quad \psi\in K.$$ 
Writing $\phi_i\in H$ for $V^*\circ \psi_i\times 1_K$ and computing using the 
multiplicativity of $V$, we find   
$$\phi_1^*\times\phi_2^*\circ S\circ \phi_3\times\phi_4=
((\psi_1\times\psi_2)^*\circ R\circ \psi_3\times\psi_4)\times 1_K.$$ 
Thus the multiplicative unitary associated with $H$ is as asserted. 
\smallskip 

  An analogous computation shows that we may replace $G_{V^*}$ by $G_R$ 
in Proposition~4.3.\smallskip 

  We now ask how the multiplicative unitary depends on the choice of 
ambidextrous Hilbert space $H\subset (\rho,\rho^2)$. Let $H'$ be another 
such Hilbert space then there is a unitary $U\in(\rho^2,\rho^2)$ such that 
$U\circ H=H'$. Let $V$ and $V'$ denote the multiplicative unitary 
operators associated with $H$ and $H'$ and set $R:=V\theta$, $R':=V'\theta'$, 
where $\theta$ and $\theta'$ are the flips on $H^2$ and $H^{'2}$, respectively. 
We compute the relation between $R$ and $R'$. Let $\psi\in H$, then 
$$R'\circ(U\circ\psi)\times 1_\rho=1_\rho\times(U\circ \psi)= 
1_\rho\times U\circ R\circ\psi\times 1_\rho.$$ 
Thus $R'=1_\rho\times U\circ R\circ U^*\times 1_\rho$. This is in contrast 
to the transformation law of $\theta$, namely 
$\theta'=u_2\circ \theta\circ {u_2}^*$ 
where $u_2:= U\times 1_H\circ 1_H\times U$ is the tensor power of $U$. 
It should be remembered however that $R'$ is intrinsically determined by 
$H'$ whereas $U$ is not. The interesting question is whether the 
associated multiplicative unitaries $V$ and $V'$ are necessarily 
equivalent and for which unitaries $U$, $U\circ H$ is ambidextrous.\smallskip 

  We present an example. Let ${\cal K}$ denote the $W^*$--tensor category 
of tensor powers of a Hilbert space $K$. Let $V\in(K^2,K^2)$ be a 
multiplicative unitary, then, as we have seen, $H:=V^*\circ K$ is an 
ambidextrous Hilbert space in $(K,K^2)$ whose associated multiplicative 
unitary is equivalent to $V$. Since we are free to choose any 
multiplicative unitary $V$, this makes it clear that the associated 
multiplicative unitaries can depend in an essential way on the 
ambidextrous Hilbert space.\smallskip 

  We may sum up the results to date in this section as follows.\smallskip 

\noindent
{\bf 4.4 Theorem} {\sl Let ${\cal K}$ be a category of Hilbert spaces with 
objects $K_n$, $n\in{\Bbb N}_0$ and $(K_0,K_0)={\Bbb C}$ equipped with 
commuting shifts $F$ and $G$ such that $GG=\hat FG$. Let $V\in(K_2,K_2)$ 
be the multiplicative unitary such that $G(\psi)=V\circ \hat F(\psi)$, 
$\psi\in K_1$, then $H:=V^*\circ F(K_1)$ is an ambidextrous Hilbert space. 
Let ${\cal H}$ be the resulting category of Hilbert spaces with commuting 
shifts $D$ and $E$ obtained by restricting $F$ and $\hat F$ to ${\cal H}$, 
then the shift $G^*$ on ${\cal K}$ defined by $G^*(\psi)=V^*\circ
F(\psi)$, 
$\psi\in K_1$ yields an isomorphism of ${\cal K},\, F,\, G$ with 
${\cal H},\, D,\, E$.}\smallskip  
 
\noindent  
{\bf Proof.} We know from Proposition 3.1 that $V$ is a multiplicative 
unitary and from Proposition 4.3 that $G^*(K_1)$ is an ambidextrous Hilbert 
space. The resulting category ${\cal H}$ of Hilbert spaces is thus
$G^*({\cal K})$. 
Since $G^*$ commutes with $F$ and by Proposition 3.1, $G^*G=\hat FG^*$, 
$G^*$ does yield the desired isomorphism.\smallskip 

\noindent
{\bf Remark.} The construction of this section deriving a multiplicative 
unitary from an ambidextrous Hilbert space is invariant under tensor 
$^*$--functors since the image of an ambidextrous Hilbert space of support 
one is again such a Hilbert space and a tensor $^*$--functor commutes 
with tensoring on the right and left.\smallskip 

  We have seen in Theorem 4.2 how an ambidextrous Hilbert space 
leads to a multiplicative unitary. There is a variant of this 
result which instead yields a representation of a multiplicative 
unitary. To motivate this result, we let $V$ be a multiplicative 
unitary on the tensor power of a Hilbert space $L$ and $W$ a 
representation of $V$ on a Hilbert space $M$. Then $V$ and $W$ 
are objects of the tensor $W^*$--category ${\cal R}(V)$ and, as 
we know, $K:=V^*\circ L\times 1_L$ is an ambidextrous Hilbert 
space in $(V,V^2)$. However, $W$ being a representation of $V$, 
$H_0:=W^*\circ M\times 1_L$ is a Hilbert space of support one 
in $(V,WV)$. Hence $H_n:=H_{n-1}\times 1_V\circ K$ is a 
Hilbert space of 
support one in $(V,WV^{n+1})$. Thus just as we have a category 
${\cal K}$ of Hilbert spaces associated with $K$, there is a category 
${\cal H}$ associated with $H_0$. We claim that tensoring on the left 
with $1_W$ restricts to a $^*$--functor from ${\cal K}$ to ${\cal H}$.
It suffices to show that $1_W\times K\circ H_0\subset H_1$ and, 
expressing $K$ and $H_0$ in terms of $L$ and $M$, this is again 
a consequence of $W$ being a representation of $V$. We have here 
an obvious generalization of the notion of ambidextrous Hilbert space 
involving $K$ and $H_0$.\smallskip 
  
  We use this as the basis of a definition. Let ${\cal T}$ be a 
tensor $W^*$--category and $\rho$ and $\sigma$ objects of ${\cal T}$. 
Let $K\subset(\rho,\rho^2)$ and $H_0\subset(\rho,\sigma\rho)$ 
be Hilbert spaces of support one and ${\cal K}$ and ${\cal H}$ 
the corresponding categories of Hilbert spaces with objects 
$K^n:=K^{n-1}\times 1_\rho\circ K$ and $H_n:=H_{n-1}\times 1_\rho\circ K$, 
$n\in{\Bbb N}_0$, respectively. We say that $H$ is $K$--ambidextrous 
if $1_\sigma\times$ restricts to a $^*$--functor from ${\cal K}$ 
to ${\cal H}$.\smallskip 

\noindent
{\bf Lemma 4.5} {\sl Let ${\cal T}$ be a tensor $W^*$--category and 
$K\subset(\rho,\rho^2)$ and $H\subset(\rho,\sigma\rho)$ be Hilbert 
spaces of support one then the following conditions are equivalent. 
\begin{description}
\item{a)} $1_\sigma\times K\subset(H_0,H_1)$,
\item{b)} $H$ is $K$--ambidextrous,
\item{c)} $H_n=1_\sigma\times K^n\circ H_0$, 
\item{d)} $H_1=1_\sigma\times K\circ H_0$.
\end{description}} 

\noindent
{\bf Proof.} ${\cal K}$ is the smallest $W^*$--subcategory containing 
$K$ and invariant under $\times 1_\rho$. Furthermore $1_\sigma\times$ 
and $\times 1_\rho$ commute, thus a) implies b). Since both sides of 
c) are Hilbert spaces of support one, b) implies c). 
c) implies d), trivially and d) implies a).\smallskip 

\noindent
{\bf 4.6 Theorem} {\sl Let $K\subset(\rho,\rho^2)$ be an ambidextrous 
Hilbert space of support one in $(\rho,\rho^2)$ and $H$ a $K$--ambidextrous 
Hilbert space in $(\rho,\sigma\rho)$. Let $F$, $G$, $E$ and $H$ denote the 
restrictions of $\times 1_\rho$, $1_\rho\times$, $1_\sigma\times$ and 
$1_{H_0}\times$, respectively, to ${\cal K}$ then there is a unique 
$V\in(K^2,K^2)$ such that 
$$G(\psi)=V\circ\hat F(\psi),\quad \psi\in K,$$ 
and a unique $W\in(H_1,H_1)$ such that 
$$E(\psi)=W\circ H(\psi),\quad \psi\in K.$$ 
$V$ is a multiplicative unitary and $W$ is a representation 
of $V$ on $H_0$.}\smallskip 

\noindent
{\bf Proof.} In view of Theorem 4.2, we need only prove the assertions 
relating to $W$. $W$ is unique and is unitary because it is given by 
$$W=\sum_iE(\psi_{i,1})\circ H(\psi_{i,1})^*$$ 
where the sum is taken over an orthonormal basis of $K$. Now since 
$H_0\in(\rho,\sigma\rho)$, for $X\in(K^m, K^n)$,
$$HG(X)=\sum_iF^{n+1}(\phi_i)\circ 1_\rho\times X\circ F^{m+1}(\phi_i)^*
=1_{\sigma\rho}\times X=EG(X),$$
where $\phi_i$ is an orthonormal basis of $H_0$. Thus $HG=EG$ and 
$W$ is a representation of $V$ by the discussion following Lemma 3.4.\smallskip 

  Theorem 4.6 does not really refer to the whole tensor $W^*$--category 
${\cal T}$ but only to the full subcategory with objects $\rho^n$ and 
$\sigma\rho^n$ and the structures induced on this subcategory by 
tensoring on the left and right by $1_\rho$ and on the left by $1_\sigma$. 
Introducing ${\cal T}$ enables us to avoid spelling out the
structures.\smallskip 

  We now show how a representation of multiplicative unitaries gives rise 
to an interchange law.\smallskip 

\noindent
{\bf 4.7 Proposition} {\sl Let $F$ and $G$ be two commuting shifts on
${\cal K}$ 
with $GG=\hat FG$ and $E$ and $H$ normal $^*$--functors from ${\cal K}$ 
to ${\cal H}$ with $EG=HG$, 
then given $X\in(H_m,H_n)$ and $Y\in(p,q)$, 
$$B^{q+1}(X)\circ EG^{m+1}(Y)=EG^{n+1}(Y)\circ B^{p+1}(X),$$ 
where $B$ is the normal $^*$--functor on ${\cal H}$ associated with $E$: 
$$B(X):\sum_iEF^n(\psi_{i,1})\circ X 
\circ EF^m(\psi_{i,1})^*, \,\,X\in(H_m,H_n).$$
Thus defining $X\times'Y:=B^{q+1}(X)\circ EG^{m+1}(Y)$ gives an 
action of the tensor $W^*$--category ${\cal K}^+$ on the category 
${\cal H}$ of Hilbert spaces. If ${\cal K}$ is given the 
tensor structure determined by $F$ then there is a unique normal tensor 
$^*$--functor $G^*$ from ${\cal K}^+$ to ${\cal K}$ such that 
$G^*(\psi):=V^*\circ F(\psi)$, $\psi\in K$, where $V$ is the
multiplicative 
unitary associated with $F$ and $G$.}\smallskip 

\noindent
{\bf Proof.} To prove the interchange law, write $EG^{m+1}=H\hat F^mG$ 
and $EG^{n+1}=H\hat F^nG$ and use the interchange law between $A$ and 
$H\hat F$ discussed before Lemma 2.4. The remarks above 
show that ${\cal K}^+$ is a tensor $W^*$--category and allow us to 
check that $G^*$ is a tensor $^*$--functor using Proposition 3.1e.\smallskip

  As we shall be considering Hilbert spaces $L\subset(K^r,K^{r+g})$ in 
Section 5, it is natural to ask 
to what extent the results can be generalized to Hilbert spaces $K$ of 
support one in $(\rho^r,\rho^{r+g})$, where we suppose, of course, that 
$g\neq 0$. The initial construction can be 
easily modified. We define inductively $K^n:=K^{n-1}\times 1_{\rho^g}\circ K$, 
the norm closed linear span being understood. $K^n$ is a Hilbert space of support 
1 in $(\rho^r,\rho^{r+ng})$. We now set 
$$(K^m,K^n):=\{X\in(\rho^{r+mg},\rho^{r+ng}): X\circ K^m \subset K^n\}.$$ 
In this way we have a $W^*$--subcategory ${\cal K}$ of Hilbert spaces
whose 
objects are of the form $\rho^{r+ng}$, $n\in{\Bbb N}_0$. ${\cal K}$ is now 
invariant under tensoring on the right by $1_{\rho^g}$. Letting $F$ be the 
restriction of this functor to ${\cal K}$ we have a shift on ${\cal K}$
that can be 
regarded as tensoring on the right by $1_K$. The analogue of Lemma 4.1 
now holds if we consider tensoring on the left by $1_{\rho^g}$ and can be 
used to define the notion of ambidextrous Hilbert space. Thus we are again
led 
to a category of Hilbert spaces ${\cal K}$ equipped with two commuting 
shifts. At this point there is the essential difference: we cannot use 
the exchange law in ${\cal T}_\rho$ to get an analogue of Theorem 4.2
unless 
$r\leq g$.\smallskip 

\noindent
{\bf 4.8 Theorem} {\sl Let $K$ be an ambidextrous Hilbert space of support
one 
in $(\rho^r,\rho^{r+g})$ where $r\leq g$. Let $F$ and $G$ denote the restrictions 
of $\times 1_{\rho^g}$ and $1_{\rho^g}\times$ to ${\cal K}$ then there is
a unique 
$V\in (K^2,K^2)$ such that 
$$G(\psi)=V\circ\hat F(\psi),\quad \psi\in K.$$ 
$V$ is a multiplicative unitary.}\smallskip 

  The proof follows that of Theorem 4.2, bearing in mind that $r\leq g$. 
If $r>g$, there is the possibility of starting with the ambidextrous 
Hilbert space $K^n$, for $n$ sufficiently large.\smallskip 

  There is another interesting general result involving Hilbert spaces in 
${\cal T}_\rho$.\smallskip 

\noindent
{\bf 4.9 Theorem} {\sl Let $K$ be a Hilbert space of support one in 
$(\rho^r,\rho^{r+g})$ and define the associated category of Hilbert 
spaces ${\cal K}$ as above. Suppose that 
$$(\rho^r,\rho^r)\times 1_{\rho^{rg}}\subset (K^r,K^r).$$ 
Then 
$$(\rho^{r+mg},\rho^{r+ng})\times 1_{\rho^{rg}}\subset (K^{r+m},K^{r+n}).$$ 
Hence defining, for $X\in(\rho^{r+mg},\rho^{r+ng})$, 
$$F(X)\psi:=X\times 1_{\rho^{rg}}\circ \psi,\quad \psi\in K^{r+m},$$ 
$F$ is a faithful $^*$--functor from the full subcategory of 
${\cal T}_\rho$ whose objects are $\rho^{r+ng}$, $n\in{\Bbb N}_0$ to 
the category ${\cal K}$ of Hilbert spaces and 
$$F(X\times 1_{\rho^{rg}})=F(X)\times 1_K.$$}

\noindent
{\bf Proof.} $K\subset (\rho^r,\rho^{r+g})$ implies 
$K^n\subset (\rho^r,\rho^{r+ng})$ and 
$K^{*n}\circ(\rho^{r+mg},\rho^{r+ng})\circ K^m\subset(\rho^r,\rho^r)$. Since 
$(\rho^r,\rho^r)\times 1_{\rho^{rg}}\subset (K^r,K^r)$, tensor the above on the 
right with $1_{\rho^{rg}}$ and compose on the right with $K^r$ and on the left 
with $K^{*r}$ to conclude that 
$$K^{*r+n}\circ(\rho^{r+mg},\rho^{r+ng})\times 1_{\rho^{rg}}\circ K^{r+m}\subset 
K^{*r}\circ(K^r,K^r)\circ K^r\subset{\Bbb C}1_{\rho^r}.$$ 
Since $K^{r+m}$ has support one, we conclude that 
$$(\rho^{r+mg},\rho^{r+ng})\times 1_{\rho^{rg}}\subset
(K^{r+m},K^{r+n}).$$
The remaining assertions are now obvious.\smallskip

   We show that the hypotheses of Theorem 4.9 with $r=g=1$ are fulfilled, 
when $V$ is a multiplicative unitary, considered as an object of ${\cal 
R}(V)$. 
In fact, $V\in(V^2,\iota(V)V)$, so $V^*\circ({\Bbb C},\iota(V))\times 1_V$ 
is a Hilbert space $H$ of support one in $(V,V^2)$.
Using the condition for $T\in(K,K)$ to be an arrow 
of $(V,V)$, we get
$$T\times 1_V\circ V^*=V^*\circ\iota(T)\times 1_V,$$ 
showing that $(V,V)\times 1_V\circ H\subset H$, as required. 
The hypothesis of Theorem 4.3 does imply for the tensor unit $\iota$ 
that $(\iota,\iota)={\Bbb C}$.\smallskip 

  Theorem 4.9 raises some interesting questions. Let ${\cal T}_{r,g}$ 
denote the full subcategory of ${\cal T}_\rho$ whose objects are of the 
form $\rho^{r+ng}$ with $n\in{\Bbb N}_0$. Then we have $^*$--functors 
$X\mapsto 1_{\rho^s}\times X\times 1_{\rho^{g-s}}$, $0\leq s\leq g$, 
defined on ${\cal T}_{r,g}$ through the ambient tensor category
${\cal T}_\rho$. 
Restricting the domain of the functors to the subcategory ${\cal K}$ and 
composing with the functor $F$ of Theorem 4.9 gives us $^*$--endofunctors 
of ${\cal K}$ taking $K^m$ to $K^{r+m}$. This raises the question of 
whether the composition with $F$ is necessary, or, more precisely, 
whether $1_{\rho^s}\times K\times 1_{\rho^{g-s}}\subset(K,K^2)$? 
In fact, we know that this is true by construction for $s=0$ and the 
basis of the definition of ambidextrous for $s=g$. When it is valid 
for some $s>r$ there is some analogue of a multiplicative unitary.\smallskip 
 
  We have seen how two commuting shifts $F$ and $G$ on ${\cal K}$ cannot 
be interpreted as tensoring on the right and tensoring on the left with 
the object $1$ unless $G=\hat F$ because the interchange law would fail 
to hold. On the other hand, we have, in this section, been using the 
interchange law in a tensor category to produce commuting shifts 
tied to multiplicative unitaries. We want, now, to show how this process 
can be reversed. We first describe a tensor category in terms of 
tensoring on the right and tensoring on the left. Let ${\cal T}$ be 
a category and give for each object $\rho$ endofunctors $F_\rho$ and 
$G_\rho$ such that for $X\in(\mu,\nu)$ and $Y\in(\pi,\rho)$, 
$$F_\rho(X)\circ G_\mu(Y)=G_\nu(Y)\circ F_\pi(X),$$ 
where this identity defines $X\times Y$ and expresses the interchange law. 
It is understood to imply that $F_\rho(\mu)=G_\mu(\rho)$ for each pair 
$\rho$, $\mu$ of objects. The set of these endofunctors is supposed 
to commute pairwise and $F_{F_\rho(\sigma)}=F_\rho F_\sigma$ and 
$G_{G_\mu(\nu)}=G_\mu G_\nu$ implying that $\times$ is associative. 
If we further require that for some object $\iota$ $F_\iota$ and $G_\iota$ 
are the identity functors then ${\cal T}$ becomes a (strict) monoidal 
category with monoidal unit $\iota$. A functor $J$ between two such 
monoidal categories is a (strict) monoidal functor if for each object $\rho$ 
of ${\cal T}$,  
$$JF_\rho=F_{J(\rho)}J,\quad JG_\rho=G_{J(\rho)}J.$$ 
To have corresponding statements for tensor categories or tensor 
$C^*$--categories or tensor $W^*$--categories we need only add the 
obvious conditions that the functors involved be linear, $^*$--preserving 
or normal as the case may be. 
\smallskip 

\noindent
{\bf 4.10 Proposition} {\sl Let $F$ and $G$ be two commuting shifts with
$GG=\hat FG$, 
then given $X\in(m,n)$ and $Y\in(p,q)$, 
$$F^{q+1}(X)\circ G^{m+1}(Y)=G^{n+1}(Y)\circ F^{p+1}(X).$$ 
Thus defining $X\times'Y:=F^{q+1}(X)\circ G^{m+1}(Y)$ gives 
a tensor $W^*$--category ${\cal K}^+$ after adjoining an irreducible
tensor unit 
${\Bbb C}$ with no arrows to any other object. If ${\cal K}$ is given the 
tensor structure determined by $F$ then there is a unique normal tensor 
$^*$--functor $G^*$ from ${\cal K}^+$ to ${\cal K}$ such that 
$G^*(\psi):=V^*\circ F(\psi)$, $\psi\in K$, where $V$ is the
multiplicative 
unitary associated with $F$ and $G$.}\smallskip 

\noindent
{\bf Proof.} Write $G^{m+1}=\hat F^mG$ and $G^{n+1}=\hat F^nG$ and 
use the interchange law between $F$ and $\hat F$. The remarks above 
show that ${\cal K}^+$ is a tensor $W^*$--category and allow us to 
check that $G^*$ is a tensor $^*$--functor using Proposition 3.1e.\smallskip
 
  Note that $\times'$ is not addition on the objects of ${\cal K}$. 
Instead we have $m\times n=m+n+1$. However as we have adjoined a 
tensor unit ${\Bbb C}$ to give ${\cal K}^+$, it is natural to renumber the 
objects by adding one and $\times'$ is then addition on the objects of 
${\cal K}^+$.\smallskip 

  We now generalize some of our results so that we can work 
with tensor $C^*$--categories rather than just tensor $W^*$--categories. 
We begin with an object $\rho$ in a strict tensor $C^*$--category and 
a Hilbert space $K$ in $(\rho,\rho^2)$ such that $K\times 1_{\rho^m}$ has 
left annihilator zero for $m\in{\Bbb N}_0$. We can define $K^n$ and 
$(K^m,K^n)$ as at the beginnning of this section to get a $C^*$--category 
${\cal K}$. ${\cal K}$ will now not be a category of Hilbert spaces but
just 
some subcategory. ${\cal K}$ is obviously an invariant subcategory for 
$\times 1_\rho$. The category ${\cal K}$ can be completed in an obvious
way 
to give a category of Hilbert spaces $\tilde{\cal K}$, say, by
identifying 
$X\in(K^m,K^n)$ with the corresponding linear map $\phi\mapsto X\circ\phi$, 
$\phi\in K^m$. The concept of ambidextrous Hilbert space is dealt with in 
the following lemma.\smallskip 

\noindent
{\bf 4.11 Lemma} {\sl Let $K\subset (\rho,\rho^2)$ be a Hilbert space 
such that $K\times 1_{\rho^m}$ has left annihilator zero for $m\in{\Bbb 
N}_0$. 
Let ${\cal K}$ be the $C^*$--category defined above, then the following 
conditions are equivalent. 
\begin{description}
\item{a)} ${\cal K}$ is an invariant subcategory for $1_\rho\times$. 
\item{b)} $1_\rho\times K\subset (K,K^2)$. 
\item{c)} $K^n=1_{\rho^{n-1}}\times K\circ K^{n-1}$, $n\in{N}_0$. 
\item{d)} $K^2=1_\rho\times K\circ K$.
\end{description}} 

\noindent
{\bf Proof.} If b) holds, then $1_\rho\times\psi\circ\psi^{'*}\times
1_\rho$ 
maps $K^2$ into $K^2$ if $\psi,\,\psi'\in K$. Hence 
$\sum_i1_\rho\times\psi_i\circ\psi_i^*\times 1_\rho$ converges in, say, 
the $s$--topology to a unitary arrow $U$ in $\tilde{\cal K}$, where 
the sum is taken over an orthonormal basis in $K$. But 
$U\circ\psi\times 1_\rho\circ\psi'=1_\rho\times\psi\circ\psi'$ and 
this proves d). c) follows from d) by induction since 
$$K\times 1_{\rho^{n-1}}\circ 1_{\rho^{n-2}}\times K
=1_{\rho^{n-1}}\times K\circ K\times 1_{\rho^{n-2}},$$
Composing on the right with $K^{n-2}$ and using the induction hypothesis 
we obtain c). But just as ${\cal K}$ is invariant under $\times 1_\rho$, 
c) implies that it is invariant under $1_\rho\times$. But b) implies 
a), trivially.\smallskip 

  We can now prove an analogue of Theorem 4.2.\smallskip 

\noindent
{\bf 4.12 Theorem} {\sl Let $\rho$ be an object in a tensor
$C^*$--category 
and $K$ an ambidextrous Hilbert space in $(\rho,\rho^2)$ such that 
$K\times 1_{\rho^m}$ has left annihilator zero for $m\in{\Bbb N}_0$. 
Let ${\cal K}$ be as above and $\tilde{\cal K}$ its completion to a 
category of Hilbert spaces, then the endofunctors on ${\cal K}$ 
determined by $\times 1_\rho$ and $1_\rho\times$ extend uniquely 
to commuting shifts $F$ and $G$ on $\tilde{\cal K}$ with $G^2=\hat FG$ and
there is 
a multiplicative unitary $V\in(\tilde K_2,\tilde K_2)$, such that 
$G(\psi)=V\circ\hat F(\psi)$, $\psi\in K$.}\smallskip 
 
\noindent
{\bf Proof.} The functor $\times 1_\rho$ defines each $K^n$ as a 
tensor power of $K$ and hence there is a unique shift $F$ on 
$\tilde{\cal K}$ such that $F(\psi_n)=\psi_n\times 1_\rho$, 
for $\psi_n\in K_n$ and $n\in{\Bbb N}_0$. If $X\in(K^m,K^n)$, then 
$$X\times 1_\rho\circ F(\psi_m)\circ\psi=
X\times 1_\rho\circ\psi_m\times 1_\rho\circ\psi=F(X\circ\psi_m)\psi.$$ 
Thus $F(X)=X\times 1_\rho$. In the same way, in view of Lemma 4.11, 
there is a unique shift $G$ on $\tilde{\cal K}$ with $G(X)=1_\rho\times
X$. 
$F$ and $G$ obviously commute and, as in Theorem 4.2, we see 
that $GG=\hat FG$. Thus there is a multiplicative unitary $V$ with the 
properties claimed.\smallskip 

  The asymmetry in the formulation of Theorem 4.12 is only apparent: 
its hypotheses imply that $1_{\rho^m}\times K$ has left annihilator 
zero for $m\in{\Bbb N}_0$.\smallskip

\section{Algebraic Endomorphisms of the Cuntz Algebra} 

  As a preliminary to our main duality result, we present in this section
results on Hilbert spaces in the Cuntz algebra and endomorphisms of that 
algebra that are `algebraic' with respect to the natural grading of the 
Cuntz algebra.\smallskip 

  When $K$ is a finite dimensional Hilbert space, ${\cal O}_K$ will denote 
the Cuntz algebra, a simple unital $C^*$--algebra introduced by Cuntz\cite{Cu}. 
When $K$ is infinite dimensional, it denotes the extended Cuntz algebra, 
a simple $C^*$--algebra introduced in \cite{CDPR}. These algebras are special 
cases of a more general construction needed in \S 6, where the Hilbert space 
is replaced by an object in a tensor $C^*$--category. ${\cal O}_K$ has a 
${\Bbb Z}$--grading derived from the automorphic action $\alpha$ of
${\Bbb T}$ 
with $\alpha_\lambda(\psi)=\lambda\psi$, $\lambda\in{\Bbb T}$. Thus the 
part of ${\cal O}_K$ of grade $k$ is given by 
$${\cal O}_K^k:=\{X\in{\cal O}_K:\alpha_\lambda(X)=\lambda^kX,\,\,
\lambda\in{\Bbb T}\}.$$

  It is this grading in the case of the Cuntz algebra ${\cal O}_K$ which 
will allow us to refer to certain Hilbert spaces in this algebra and 
endomorphisms of this algebra as being `algebraic'. Here $K$ can be a finite 
dimensional or infinite dimensional Hilbert space and in the latter case
${\cal O}_K$ is the extended Cuntz algebra introduced in \cite{CDPR}. 
The results can be roughly summarized by saying 
that computations involving these algebraic objects reduce to problems involving 
linear operators between Hilbert spaces that can be identified a priori. 
These results will prove useful in other contexts when dealing with 
concrete endomorphisms on the Cuntz algebra.\smallskip 

  The problem has its origins in the relation between the Cuntz algebra 
${\cal O}_K$ and the tensor $W^*$--category of bounded linear mappings
between 
tensor powers of $K$. The Cuntz algebra is obtained from the category by 
factoring out the operation of tensoring on the right by $1_K$ whilst 
the operation of tensoring on the left is retained in the shape of 
the canonical endomorphism $\rho_K$. Then a direct sum is taken over 
the grading and finally the algebra is completed in the unique $C^*$--norm. 
This raises certain questions when proving results using the Cuntz algebra. 
A result may be purely algebraic in nature involving only the algebraic 
part of the Cuntz algebra. The manipulations may have been simplified by 
factoring out the operation of tensoring on the right by $1_K$. At the 
same time their significance may have been obscured by the asymmetric 
treatment of tensoring on the two sides.\smallskip

  The problems treated in this section illustrate both the analytic 
and the algebraic aspects of the problem and we have had more success 
with the analytic aspects proving that the set of solutions of these 
problems involves only the algebraic part.\smallskip  

   If $H$ is a Hilbert space in a $C^*$--algebra ${\cal A}$,
$Y\in{\cal A}$, and $L^1(H)$ denotes the trace class operators on $H$, 
then there is a unique continuous linear mapping $T\mapsto \text{Tr}_H(YT)$ 
from $L^1(H)$ to ${\cal A}$ such that 
$$\text{Tr}_H(Y\psi'\psi^*)=\psi^*Y\psi',\quad \psi,\psi'\in M.$$ 
The norm of this mapping is $\leq\|Y\|$ and Tr$_H(YT)=$Tr$_H(TY)$. As the 
notation suggests, we are taking a partial trace relative to $H$. 
\smallskip 

\noindent
{\bf 5.1 Lemma.} {\sl Let ${\cal C}$ be a ${\mathbb Z}$--graded
$^*$--subalgebra of 
${\cal O}_K$, i.e. ${\cal C}$ is generated by the subspaces ${\cal
C}^k:={\cal C}\cap{\O_K}^k\ .$ 
Then ${\cal C}^k\subset H$ for some $k$ and some Hilbert space $H$ in
${\cal O}_K$ 
of dimension $>1$ implies ${\cal C}^{nk}=0$ for $n\in{\mathbb Z}\neq
0$.}\smallskip

\noindent{\bf Proof.} If $n\in{\mathbb N}$, then ${\cal C}^{nk}\subset
H^n$ so it \
suffices to show ${\cal C}^k=0$. Let $\psi\in{\cal C}^k$, then 
$\psi\psi\psi^*\in{\cal C}^k\subset H$. Hence
$\psi^*\psi\psi\psi^*\in{\mathbb C}I$, 
but $\psi\psi^*$ is a multiple of a minimal projection in $(H,H)$ so 
$\psi=0$.\smallskip

   In our applications, the ${\mathbb Z}$--grading will be that introduced
above 
and $H=K^k$. We next give results on computing the relative commutant of 
certain Hilbert spaces of support $I$ in the Cuntz algebra. This amounts to 
the same thing as determining the fixed points under the inner endomorphism 
generated by the Hilbert spaces in question. If $\lambda$ is an endomorphism 
of the Cuntz algebra, then $(\lambda,\lambda)$ is just the relative commutant 
of $\lambda(K)$ and we give our result in a generality to include computing 
certain spaces of intertwining operators. We consider Hilbert spaces of support 
$I$ that are {\it algebraic} in the sense that they are contained in some 
$(K^r,K^{r+g})$. $g$ will be referred to as the grade of the Hilbert space 
and we shall only consider the case that $g\geq 1$. Indeed if $K$ is 
finite--dimensional $g\geq 0$ and $g=0$ implies that the Hilbert space has 
dimension one.  The minimal value for $r$ will be referred to as 
the rank of the Hilbert space. If $K$ is finite dimensional, then every 
such Hilbert space is of the form $WK$, where $W\in (K^{r+1},K^{r+g})$ is 
an isometry. In infinite dimensions, we need to consider coisometries $W$, too. 
As these Hilbert spaces have a fixed grade, 
the endomorphisms they generate commute with $\alpha_\lambda$ and their 
fixed--point algebras are graded $C^*$--subalgebras of ${\cal O}_K$.
The basic observation is that if $L$ and $M$ are such Hilbert spaces of grade 
$g$ and rank $q$ and $r$, respectively, then 
$$L^*(K^m,K^{m+k})M\subset (K^{m-1},K^{m-1+k}),\quad m\geq q-k+g,r+g,
\footnote{For convenience, we regard expressions involving the compositions of 
linear spaces as referring to the norm closed linear span.}$$ 
$$L^*(K^m,K^{m+k})M\subset (K^r,K^{r+k}),\quad r+g\geq m,\,q-r\leq k,$$
$$L^*(K^m,K^{m+k})M\subset (K^{q-k},K^q),\quad q-k+g\geq m,\,q-r\geq k.$$ 
This will be used in the following way. Pick $\varphi\in L$ 
and $\psi\in M$ of norm $\leq 1$ and consider the linear mapping 
$X\mapsto \varphi^*X\psi$ of norm $\leq 1$ on ${\cal O}_K$. Let $\Phi$ 
denote a limit point of iterates of such mappings in the pointwise 
weak operator topology of some locally normal representation of ${\cal O}_K$ 
on some Hilbert space ${\cal H}$. A priori $\Phi$ maps ${\cal O}_K$ into 
${\cal B}({\cal H})$, however since the subspaces of the form $(K^m,K^n)$ 
are weak operator closed in such representations\cite{CDPR}, $\Phi$ 
will map ${\cal O}_K^k$ into $(K^r,K^{r+k})$ or $(K^{q-k},K^q)$ according 
as $k\geq(q-r)$ or $k\leq (q-r)$. Hence $\Phi$ maps ${\cal O}_K$ into 
itself.\smallskip 

   We now consider the following intertwining problem. Given a bounded 
linear mapping $Y\in (M,L)$ of norm $\geq 1$, find the set ${\cal C}$ of 
elements $X$ of ${\cal O}_K$ satisfying one of the following equivalent conditions 
\begin{description}  
\item{a)} $X\psi=Y\psi X,\quad \psi\in M$,
\item{b)} $\psi^{'*}X\psi=\psi^{'*}Y\psi X,\quad \psi\in M,\,\,\psi'\in L$,
\item{c)} $\psi^{'*}X=X\psi^{'*}Y, \quad \psi'\in L$, 
\item{d)} $X^*\psi'=Y^*\psi'X^*, \quad \psi'\in L$, 
\item{e)} $X=Y\rho_M(X)$.
\end{description}
Notice that as we have chosen Hilbert spaces $L$ and $M$ of equal grade, 
${\cal C}$ is stable under the automorphisms $\alpha_\lambda$ defining 
the ${\mathbb Z}$--grading. To compute ${\cal C}$ it therefore suffices to 
compute ${\cal C}^k$ for each $k$. The first step is to use an appropriate 
mapping $\Phi$. If we can pick $\varphi\in L$ and $\psi\in M$ of norm 
$\leq 1$ such that $\varphi^*Y\psi=I$, 
then by b), the map $X\mapsto\varphi^*X\psi$ leaves ${\cal C}$ 
pointwise invariant, and letting $\Phi$ be a limit point of iterates of 
this particular mapping, we conclude that 
$${\cal C}^k\subset (K^r,K^{r+k}),\quad k\geq (q-r),$$
$${\cal C}^k\subset (K^{q-k},K^q),\quad k\leq (q-r).$$
In general, we could define $\Phi$ to be a pointwise weak operator limit 
point of $X\mapsto \varphi_n^{*n}X\psi_n^n$, where the $\varphi_n$ and $\psi_n$ 
are chosen of norm $\leq 1$ such that $(\varphi_n^*Y\psi_n)^n\to I$ as $n\to\infty$, 
this being possible since $\|Y\|\geq 1$ by assumption. 
But we want to go further and reduce the problem of computing ${\cal C}^k$ to 
a purely local problem.\smallskip 

\noindent
{\bf 5.2 Proposition} {\sl Let $L$ and $M$ be algebraic Hilbert spaces 
of equal grade and rank $q$ and $r$, respectively and $Y\in(M,L)$ of norm 
$\geq 1$. 
Let ${\cal C}$ denote the set of $X\in {\cal O}_K$ such that 
$$X\psi=Y\psi X,\quad \psi\in M,$$ 
then if $k\geq q-r$, $X\in {\cal C}^k$ if and only if $X\in (K^r,K^{r+k})$ 
and one of the following equivalent conditions hold 
\begin{description}
\item{a)} $X\vartheta(K^r,M)=Y\vartheta(K^{r+k},M)X$. 
\item{b$_n$)} $X\text{Tr}_{M^n}(TY^{\times n}\vartheta(K^r,M))=\text{Tr}_{L^n}
(Y^{\times n}TY\vartheta(K^{r+k},M))X$, $T\in~L^1(L^n,M^n).$ 
\end{description}
Here $L^1$ is used to denote the set of trace class operators. If $k\leq q-r$, 
then $X\in (K^{q-k},K^q)$ and we need only replace $r$ by $q-k$ in the above. 
}\smallskip

\noindent
{\bf Proof.} We have already seen that $X\in (K^r,K^{r+k})$ if $k\geq q-r$ 
and a) now follows noting that: 
$$X\vartheta(K^r,M)=Y\rho_M(X)\vartheta(K^r,M)=Y\vartheta(K^{r+k},M)X.$$ 
Conversely, a) implies $X=Y\rho_M(X)$ since $X\in (K^r,K^{r+k})$. 
Now a) also implies that 
$$X\psi^{'*}Y\vartheta(K^r,M)\psi=\psi^{'*}Y\vartheta(K^{r+k},M)Y\psi X,\quad \psi'\in L,
\,\,\psi\in M.$$ 
This is b$_1$) for rank one operators $T$ and hence equivalent to b$_1$). 
On the other hand, a) follows from b$_1$) since $M$ and $L$ have support one. 
The same argument shows that b$_n$) is equivalent to b$_{n-1}$), completing 
the proof.\smallskip    
 
   Let us comment on these conditions: a) is 
a simple canonical condition that 
already serves to make the basic point that ${\cal C}^k$ is determined by 
intertwining conditions between fixed tensor powers of $K$ and is in this 
sense algebraic. However the permutation operators map between higher 
tensor powers of $K$ than is really necessary if $X\in (K^r,K^{k+r})$ is 
to intertwine. By using the partial trace, we can reduce the powers of 
the tensor spaces involved at the cost of increasing the number of 
intertwining relations. In fact, for $n$ sufficiently large, the operators 
involved on the left hand side are in $(K^r,K^r)$ and those on the right 
hand side in $(K^q,K^q)$.\smallskip
 
  In concrete cases, the following strategy for computing ${\cal C}^k$ 
proves useful. Let $X\in{\cal C}^k$ and
$V\in(K,M)$; in practice, $V$ can usually be picked unitary. Then 
$$\vartheta(K^n,M)\rho^n(V)=V\vartheta(K^n,K),\quad n\in{\Bbb N}_0,$$ 
where we have written $\rho$ for $\rho_K$. Since $X\in(K^r, K^{r+k})$, 
$X\rho^r(V)=\rho^{r+k}(V)X$, and using a) of Proposition 5.2, we get 
$$XV\vartheta(K^r,K)=YV\vartheta(K^{r+k},K)X.$$ 
If $V$ has a right inverse, we can, conversely, deduce a) of Proposition 5.2 
from this equation. If wished, the permutations operators can be 
eliminated in favour of the endomorphism $\rho$. In fact, since 
$\vartheta(K^{r+k},K)X=\rho(X)\vartheta(K^r,K)$, we get 
$$XV=YV\rho(X),$$ 
but this is best derived directly from d), above. Similarly, if 
$U\in(K,L)$, we may conclude that 
$$UX=\rho(X)UY.$$
This is equivalent to a) of Proposition 5.2, if $V$ has a left inverse.\smallskip 
 
  After these results, let us try and clarify whether more might be expected  
by relating the set of solutions to questions posed entirely in terms 
of a category of Hilbert spaces and hence independent of the 
identifications used to define the Cuntz algebra.\smallskip 

  In place of ${\cal O}_K$ we consider a category ${\cal K}$ of Hilbert
spaces 
with objects ${\Bbb N}_0$ and equipped with a shift $F$ to be thought of
as 
the tensor powers of a Hilbert space $K$, as described in detail in
Section 2. 
Instead of considering a Hilbert spaces $L\subset (K^q,K^{q+g})$ and 
$M\subset (K^r,K^{r+g})$, we consider another such category ${\cal L}$
with 
a shift $G$ and two normal $*$--functors $J$ and $J'$ from ${\cal L}$ to 
${\cal K}$ such that $JG=F^gJ$, $J'G=F^gJ'$ and $J(0)=q$ and $J'(0)=r$. 
We consider natural transformations $t\in(J',J)\cap(\hat F,\hat F)$. 
As we know from computations in Section 2, such a natural 
transformation is uniquely determined by
$t_0:=X\in (K^q,K^r)$ satisfying 
$$F^g(X)=\sum_jJ(\psi_j)\circ X\circ J'(\psi_j)^*,$$
where the sum is taken over an orthonormal basis of $M$. If we 
write this in the Cuntz algebra we obtain our condition e), 
$X=Y\rho_M(X)$, where $Y=\sum_jJ(\psi_j)J'(\psi_j)^*$ in the 
Cuntz algebra and is unitary. There is no difficulty in 
generalizing to include cases where $Y$ is not unitary.\smallskip 

  We learn from this that it is quite natural to expect solutions 
of grade $r-q$. Furthermore, by composing $J$ or $J'$ with tensoring on the 
right by $1_K$, we can replace $q$ by $q+1$ or $r$ by $r+1$.  
Thus we have potential solutions for any grade. 
We see, therefore that the identifications involved 
in defining the Cuntz algebra mean that one problem at the level of the 
Cuntz algebra involves a countable set of problems at the level of ${\cal 
T}_K$.\smallskip 
 
   We conclude that the results obtained using the spaces $(K^m,K^n)$ 
in the Cuntz algebra are the best that can be expected in complete  
generality. However, we now show how the 
estimates on the localization of ${\cal C}^k$ can be improved under 
conditions involving the relative localization of $YM$ and $K^g$ or 
$Y^*L$ and $K^g$. Note that $M^{n*}K^{gn}\subset (K^r,K^r)$ for all 
$n$.\smallskip 

\noindent
{\bf 5.3 Lemma} {\sl Let $m$ denote the smallest integer $\geq\frac{r}{g}$ 
and $\ell$ the smallest integer $\geq{\frac{q}{g}}$. If the weak operator 
closed linear span 
of the $(L^*Y)^nK^{gn}$, $n\geq m$, in $(K^r,K^r)$ contains 
$I$ then ${\cal C}^k\subset K^k$ for $k\geq q$ and ${\cal C}^k\subset (K^{q-k},K^q)$ 
for $q\geq k\geq q-r$. Similarly, if the weak operator closed linear 
span of the $(M^*Y^*)^nK^{gn}$, $n\geq\ell$, in $(K^q,K^q)$ 
contains $I$ then ${\cal C}^k\subset K^{*-k}$ for $k\leq -r$ and 
${\cal C}^k\subset (K^r,K^{r+k})$ for $-r\leq k\leq q-r$.}\smallskip

\noindent
{\bf Proof.} We know that ${\cal C}^k\subset (K^r,K^{r+k})$ if $k\geq q-r$. 
But $L^n$ has rank $q$ and grade $ng$ and $K^{ng}$ has rank $0$ and grade 
$ng$. Thus from previous computations
$$L^{n*}{\cal C}^kK^{ng}\subset K^k,\quad k\geq q,\,\, ng\geq r,$$ 
$$L^{n*}{\cal C}^kK^{ng}\subset (K^{q-k},K^q),\quad q\geq k\geq q-r,\,\,
ng\geq r.$$ 
But if $X\in {\cal C}$,  
$$L^{n*}XK^{gn}=X(L^*Y)^nK^{gn}.$$ 
Thus if the weak operator closed linear span of the $(L^*Y)^nK^{gn}$ 
contains $I$, $X$ will be in the weak operator closed linear span of the 
$L^{n*}XK^{gn}$ and the first part follows. The second part can be proved 
similarly or deduced from the first by using $X^*$ and $Y^*$ in place 
of $X$ and $Y$.\smallskip 

   Recalling Lemma 5.1 at this point, we get the following corollary.\smallskip 

\noindent
{\bf 5.4 Corollary} {\sl Suppose $L=M$ and $Y$ is a projection or a unitary. 
Let $m$ denote the smallest integer $\geq\frac{q}{g}$ and suppose the 
weak operator closed linear span 
of the $(L^*Y)^nK^{gn}$, $n\geq m$, in $(K^r,K^r)$ contains $I$,  
then ${\cal C}^k=0$ for $k\geq q$ and $k\leq -q$ and 
${\cal C}^k\subset(K^{q-k},K^q)$ for $q\geq k\geq0$.}\smallskip 

\noindent
{\bf Proof.} When $Y$ is a projection, we need only remark that ${\cal C}$ 
is a ${\mathbb Z}$--graded $^*$--subalgebra of ${\cal O}_K$. If $Y$ is
unitary, then $X\in{\cal C}^k$ implies $XXX^*\in{\cal C}^k$ and this 
is all that is used in Lemma 5.1.\smallskip 

   To give a simple example, $\rho^r(K)$ is an algebraic Hilbert space of grade 
one and rank $r$. Its relative commutant is $(\rho^r,\rho^r)$ which was 
shown in \cite{DR2}, using techniques similar to those above, to be equal 
to $(K^r,K^r)$. In this case, $\rho^r(K^{n*})K^{n}$ is the space of finite rank 
operators on $K^r$ for $n\geq r$ and the space of compact operators from 
$\rho^{r-n}(K^n)$ to $K^n$ if $1\leq n\leq r$.\smallskip 

  The theory of multiplicative unitaries provides us with further examples 
of algebraic Hilbert spaces $L$ and $M$ of equal grade, where the weak 
operator closure of $L^*M$ and hence of $M^*L$ contains $I$. For example, 
if $V\in(K^2,K^2)$ is a regular multiplicative unitary then the weak 
operator closures of $K^*VK$, $K^*\vartheta V\vartheta K$ and $K^*V\vartheta K$  
in $(K,K)$ are even $^*$--algebras containing the unit \cite{BS}. 
Of course, $K$ could be replaced here by any other algebraic 
Hilbert space of support $I$. Note that if $L_i$ and $M_i$ are algebraic 
Hilbert spaces such that $I$ is in the weak operator closure of 
$L^*_iM_i$, $i=1,2$, then $I$ is also in the weak operator closure of 
$L^*_1L^*_2M_2M_1$.  

   Let us now return to the special case used to motivate our 
basic intertwining relation. If we take $Y=1_M$ then ${\cal C}$ is just 
the relative commutant of $M$ and it is of interest to ask when $M$ has trivial 
relative commutant. ${\cal C}^0$ will reduce to the complex numbers by b$_1$) of 
Proposition 5.2 if $M^*\vartheta(K^r,M)M$ has trivial commutant in 
$(K^r,K^r)$. We again have examples with $r=1$ drawn from the theory 
of a regular multiplicative unitary and this leads to the following result.\smallskip 

\noindent
{\bf 5.5 Proposition} {\sl Let $V$ be a regular multiplicative unitary in 
$K^2$, then the following Hilbert spaces have trivial relative commutant 
in ${\cal O}_K$: $V^*K$ and $\vartheta V\vartheta K$. In the case of the Hilbert spaces 
$\vartheta VK$ and $V^*\vartheta K$, the relative commutants are the commutants of 
$K^*VK$ and $K^*\vartheta V^*\vartheta K$ in $(K,K)$, respectively.}
\smallskip

\noindent
{\bf Proof.} We need only remark that in each case we know that 
${\cal C}^k=0$ for $k\geq 1$. Furthermore, ${\cal C}^0$ is, 
in each case, as claimed since, for a Hilbert space of the 
form $UK$ with $U\in(K^2,K^2)$ unitary, b$_1$ of Proposition 5.2 
just reduces to saying that ${\cal C}^0$ is the commutant of the 
first component of $U\vartheta$. 
\smallskip 

   Following \cite{BS}, we denote $K^*VK$ and $K^*\vartheta V^*\vartheta K$ 
by ${\cal A}(V)$ and $\hat{\cal A}(V)$, respectively. If $V$ is a 
regular multiplicative unitary, these algebras are actually $^*$--algebras  
\cite{BS}.\smallskip

   We now come to the second application we had in mind, namely to study 
intertwiners between certain endomorphisms 
of the Cuntz algebra. We say that an endomorphism $\tau$ has grade $g$ if 
$\tau(K)$ has grade $g+1$ and is algebraic of rank $r$ if $\tau(K)$ is 
algebraic of rank $r$. There is a unique unitary $V$ such that $\tau(\psi)=
V\psi$, for $\psi\in K$ and $\tau$ has rank $r$ if and only if $r$ is the 
smallest integer such that $V\in (K^{r+1+g},K^{r+1+g})$. If $\sigma$ is an 
algebraic endomorphism of grade $f$ and rank $q$ then  the composition 
$\sigma\tau$ is of grade $f+g$ and rank $\leq q+r+g+(r+g)f$.\smallskip 

   Now suppose that we have endomorphisms $\sigma$ and $\tau$ as above 
of equal grade $g$, then the space of intertwiners $(\tau,\sigma)$ 
is ${\mathbb Z}$--graded and if $X\in (\tau,\sigma)$, then 
$$X\psi=Y\psi X,\quad \psi\in \tau(K),$$ 
where $Y\in(\tau(K),\sigma(K))$ is the unitary taking $\tau(\psi)$ to 
$\sigma(\psi)$ for each $\psi\in K$. 
Thus the analysis of Proposition 5.2 holds. In particular, we have 
$$(\tau,\sigma)^k\subset (K^r,K^{r+k}),\quad k\geq q-r,$$ 
$$(\tau,\sigma)^k\subset (K^{q+k},K^q),\quad k\leq q-r.$$

\section{An Algebraic Version of Takesaki--Tatsuuma Duality} 

After these results on algebraic Hilbert spaces and endomorphisms, our aim is 
to describe an algebraic model for a dual of a multiplicative unitary. 
We present a duality result for `locally compact' multiplicative unitaries in 
terms of the $C^*$--algebra generated by the regular representation considered 
as an object in the tensor $C^*$--category of representations of the 
multiplicative unitary. Thus we get, in particular, an algebraic version 
of a duality result for the representation categories of locally compact 
groups.\smallskip

We recall that ${\cal O}_{H}$ is a simple $C^*$--algebra and every 
unitary operator $X\in(H, H')$
extends to a unital morphism ${\cal O}_H\to{\cal O}_{H'}$, 
which we denote by $\lambda_X$.\smallskip

  Let $V$ be a regular multiplicative unitary acting on $K^2$ and $W$ 
a  representation contained in  $ M((H,H)\otimes{\cal A}(V))\ ,$ 
the multiplier algebra of the minimal tensor product 
$(H,H)\otimes{\cal A}(V)\ .$ Let us identify ${\cal A}(V)\rho_K({\cal O}_H)=
{\cal A}(V)\otimes{\cal O}_H\ .$ Then 
$$\lambda_{\vartheta_{H, K}W}(H^r, H^s){\cal A}(V)+{\cal A}(V)
\lambda_{\vartheta_{H, K}W}(H^r, H^s)\subset{\cal A}(V)\rho_K(H^r, H^s).$$ 
Here in the definition of $\lambda_{\vartheta_{H,K}W}$ we consider 
$H$ and $K$ as Hilbert spaces of support $I$ in some ${\cal B}({\cal H})$
and 
regard $\vartheta_{H,K}W$ as mapping $H$ onto $\vartheta_{H,K}WH$. 
 The arguments of \cite{Cu2} or \cite{P};
Section 6 generalize to show that  the  monomorphism
$\lambda_{\vartheta_{H, K}W}: {\cal O}_H\to M({\cal A}(V)\otimes{\cal O}_H)$ 
defines a coaction of ${\cal A}(V)$ on ${\cal O}_H\ ,$ i.e.\ 
$\lambda_{\vartheta_{H, K}W}$ is a unital $^*$--homomorphism  satisfying 
$${\delta}\otimes i\circ\lambda_{\vartheta_{H, K}W}
=i\otimes\lambda_{\vartheta_{H, K}W}\circ\lambda_{\vartheta_{H, K}W}\ ,$$ 
where ${\delta}$ is the coproduct of ${\cal A}(V)$ induced 
by the adjoint action of\linebreak $\vartheta_{K,K}V$ 
\cite{BS}. The corresponding fixed point algebra is

$${\cal O}^W=\{T\in{\cal O}_H: \lambda_{\vartheta_{H, K}W}(T)=\rho_K(T)\}\ .$$

Direct computations show that
${\cal O}^W\cap(H^r, H^s)=(W^{\times r}, W^{\times s})$. 
  We recall that any object $W$ in a tensor $C^*$--category ${\cal T}$ has 
a canonically associated $C^*$--algebra 
${\cal O}_W$ with a unital endomorphism $\rho_W$\cite{DR}. This
construction 
can be used in two quite different ways. On the one hand it provides one 
with a large class of model endomorphisms with rather well defined properties. 
In favourable cases, the associated $^*$--functor $F_W$ from ${\cal T}_W$,
the 
tensor $C^*$--category whose objects are the tensor powers of $W$ with arrows 
taken from ${\cal T}$, to the tensor $C^*$--category whose objects are the
powers 
of $\rho_W$ and whose arrows are intertwining operators is even an isomorphism.  
This illustrates the second aspect of the construction that it also encodes 
properties of the object $W$ in question. $W$ is $C^*$--amenable when $F_W$ 
is an isomorphism, in the terminology of \cite{LR}, where related notions of 
amenability are discussed.\smallskip 

  Now ${\cal T}_W$ carries an automorphic action of the circle group
${\Bbb T}$ 
defined by 
$$\alpha_\lambda(T):=\lambda^{s-r}T,\quad T\in(W^r,W^s),$$ 
and the induced automorphic action of ${\Bbb T}$ on $({\cal O}_W,\rho_W)$
is 
also denoted by $\alpha$. The spectral subspaces of the action make 
${\cal O}_W$ into a ${\Bbb Z}$--graded $C^*$--algebra:
$${\O_W}^k=\{T\in \O_W: \alpha_{\lambda}(T)=\lambda^k T\
,\quad\lambda\in{\mathbb T}\}\ .$$

 The construction is functorial so, given a $^*$--functor $F$ from ${\cal
T}$, 
it yields a morphism $F_*:{\cal O}_W\to{\cal O}_{F(W)}$ of $C^*$--algebras 
intertwining the canonical endomorphisms and the actions of ${\Bbb T}$. In 
particular, if we have a faithful functor into a tensor $C^*$--category of 
Hilbert spaces as is the case for the categories ${\cal R}(V)$ or 
${\cal C}(V)\ ,$ then it yields an inclusion 
${\cal O}_W\subset {\cal O}_{H}$ of $C^*$--algebras  such that 
$\rho_{H}\upharpoonright {\cal O}_W=\rho_W\ .$ Here $H=F(W)$ is the 
Hilbert space of $W$. If ${\cal O}_W$ has trivial relative
commutant in ${\cal O}_H$ then the group of automorphisms of 
${\cal O}_H$ leaving ${\cal O}_W$ pointwise fixed can be identified with 
$$G_W:=\{U\in{\cal U}(H):TU^{\times r}=U^{\times s}T,\quad 
T\in(W^{\times r},W^{\times s}),\,\,r,s\in{\Bbb N}\}.$$ 
Furthermore, 
$$(\rho^r_W,\rho^s_W)=(H^r,H^s)\cap{\cal O}_W.$$

Returning to the fixed point algebra under the above action, we have
${\cal O}_W\subseteq{\cal O}^W\ .$ As in the case of a group action,  
equality follows from an amenability condition on $V\ .$\smallskip

\noindent{\bf 6.1 Theorem} {\sl Let $W\in M((H,H)\otimes{\cal
A}(V))$ be a unitary representation of $V$ 
 and suppose that there is
an invariant mean $m$ on $({\cal A}(V), {\delta})\ .$ Then there is a
conditional expectation $E: {\cal O}_H\to{\cal O}^W$ satisfying $E(H^r,
H^s)=(W^{\times r}, W^{\times s})\ .$ In particular, ${\cal O}_W={\cal O}^W\
.$}\smallskip

\noindent{\bf Proof.} 
 We extend $m$  to the multiplier algebra of ${\cal A}(V)$ via strict continuity.
Let $\omega$ be a normal state of some faithful representation $\pi$ of $\cal O_H$ 
where $\pi(H)$ has support $I$. Then $i\otimes\omega$ induces a strictly 
continuous positive map from $M({\cal
A}(V)\otimes{\cal O}_H)$ to $M({\cal A}(V))\ ,$ and setting

$$\omega(E(T)):=m\circ i\otimes\omega\circ\lambda_{\vartheta_{H, K}W}(T)\ ,
\quad T\in{\cal O}_H\ ,$$
gives a positive map $E$ of norm one  from ${\cal O}_H$ to ${\cal B}({\cal
H}_\pi)$ satisfying
$E(AT)=AE(T)\ ,\quad A\in{\cal O}^W\ , T\in{\cal O}_H\ .$
Now the arguments of \cite{P}; Proposition 6.5, except that  
$\rho_K(H^s)^*\lambda_{\vartheta_{H, K}W}(H^r, H^s)\rho_K(H^r)\subseteq
M({\cal A}(V))$ here,   show that $E$ is the desired conditional
expectation.\smallskip

In particular, if $V$ is compact and $T\in K$ is a fixed normalized vector then
$m(A) I=T^*AT\ ,\  A\in{\cal A}(V)\ ,$ is the unique Haar measure on 
${\cal A}(V)\ ,$ the  conditional expectation corresponding to the
representation $W$ is $E(X)=T^*\lambda_{\vartheta_{H,K}W}(X)T$.\smallskip

\noindent
{\bf 6.2 Theorem} {\sl Let $V$ be a regular multiplicative unitary, then 
$V$ is $C^*$--amenable as an object of ${\cal R}(V)$, i.e.\ 
$$(\rho_V^r,\rho_V^s)=(V^r,V^s),\quad r,s\in{\Bbb N}_0.$$} 

\noindent
{\bf Proof.} The pentagon equation expressing the multiplicativity 
of $V$ is equivalent to $V^*K\subset(V,V^2)$. $V^*K$ has trivial relative 
commutant, by Proposition 5.5, if $V$ is regular. Thus 
${\cal O}_V'\cap{\cal O}_K={\Bbb C}$ and, consequently, 
$$(\rho_V^r,\rho_V^s)=(K^r,K^s)\cap{\cal O}_V.$$ 
On the other hand, a computation, cf.\ \S 6 of \cite{R}, shows that 
$$(V^r,V^s)=(K^r,K^s)\cap{\cal O}^V,$$ 
where ${\cal O}^V=\{X\in{\cal O}_K:\lambda_{\vartheta
V}(X)=\rho_K(X)\}$.\smallskip

   It is not clear whether $V$ is $C^*$--amenable as an object of 
${\cal C}(V)$, when $V$ is regular. The analogous proof does not work.  
The pentagon equation is now equivalent to 
$$V\rho_K(K)=V\vartheta K\subset (V,V^2).$$ 
But if $V$ is regular, the relative 
commutant of $V\vartheta K$ is, by Proposition 5.5, the commutant 
of $K^*\vartheta V\vartheta K$, the first component of $V$, in $(K,K)$. 
This difference between 
${\cal C}(V)$ and ${\cal R}(V)$ relates to the alternative definition of 
tensor product pointed out in the introduction. With the 
alternative tensor product for ${\cal C}(V)$, we would find 
$\vartheta V\vartheta K\subset (V,V^2)$ and $\vartheta V\vartheta K$ does have trivial 
relative commutant in ${\cal O}_K$.\smallskip 

   In virtue of Theorem 6.2, the model endomorphism $({\cal O}_V,\rho_V)$ is a 
natural candidate for a dual of $V$. So we  pose the following question. 
When does a pair $({\cal A},\hat\rho)$ consisting of a unital $C^*$--algebra and
 a unital endomorphism arise from a system of the form 
$({\cal O}_V, \rho_V)$?\smallskip

   At the same time we have seen that the Cuntz algebra allows a simple 
description of a large variety of interesting model endomorphisms, giving 
rise to systems of the form $({\cal O}_K,\lambda_R)$, where ${\cal O}_K$
is 
the (extended) Cuntz algebra over the Hilbert space $K$ and $\lambda_R$ the 
algebraic endomorphism determined by the unitary operator $R$, 
$R\in (K^{r+1},K^{r+1+g})$ for an algebraic endomorphism of grade $g$ and 
rank $r$. In fact, we would like our model systems to combine two features: 
they should be of the form $({\cal O}_\rho,\hat\rho)$, where $\rho$ is an 
object in a tensor $C^*$--category. This means that ${\cal O}_\rho$ is 
generated by the intertwiners between the powers of $\hat\rho$. The second 
feature is that there should be a Hilbert space $K$ of intertwiners 
between powers of $\hat\rho$.

  Now the Cuntz algebra and the extended Cuntz algebra ${\cal O}_K$ are
derived 
from the tensor $W^*$-category of Hilbert spaces whose objects are the 
tensor powers of $K$. If we begin, as in previous sections, simply with 
a $W^*$--category of Hilbert spaces ${\cal K}$ we need a shift $F$ 
to give ${\cal O}_K$. Giving a second commuting shift $G$ amounts to 
giving an algebraic endomorphism $\lambda_R$ of grade zero and rank one. 
Here $R\in(K^2,K^2)$ is determined by $G(\psi)=R\circ F(\psi)$, 
$\psi\in K$. Natural transformations yield intertwiners of endomorphisms: 
more precisely, if $t\in(G^r,G^s)\cap (\hat F^r,\hat F^s)$ then 
$t_0\in(\lambda_R^r,\lambda_R^s)$ but it is not clear whether all 
intertwiners arise in this way.\smallskip 

  Looking at Theorem 4.4 in this light gives us the following result.\smallskip 

\noindent 
{\bf 6.3 Theorem} {\sl Let $V$ be a multiplicative unitary on $K^2$, viewed 
as an element of ${\cal O}_K$. Let $H:=V^*K$ and let $\lambda^*$ be the 
induced isomorphism of ${\cal O}_K$ onto ${\cal O}_H\subset{\cal O}_K$, 
then $\rho\lambda^*=\lambda^*\lambda_R$, where $R=V\theta$ and $\rho$ 
is the canonical endomorphism of ${\cal O}_K$.}\smallskip 

\noindent
{\bf Remark} In fact, ${\cal O}_H={\cal O}_V$ since 
$(V^r,V^s)\times 1_K\subset(H^r,H^s)$ for $r,s>0$ but this still leaves 
the finer details open on how the intertwiner spaces $(V^r,V^s)$ sit in 
${\cal O}_H$.\smallskip 

  We now take up the situation of Theorem 4.9. Thus we suppose that 
${\cal T}_\rho$ is a tensor $W^*$--category whose objects are the 
powers of $\rho$ with a Hilbert space $K$ of support one in 
$(\rho^r,\rho^{r+g})$ and suppose that 
$(\rho^r,\rho^r)\times 1_{\rho^{rg}}\subset(K^r,K^r)$. We now compare 
$C^*$--algebra ${\cal O}_\rho$ with the Cuntz algebra ${\cal O}_K$. 
Obviously, in constructing ${\cal O}_K$, we use only spaces of 
arrows between objects of the form $\rho^{r+ng}$ with $n\in{\Bbb N}_0$. 
But these spaces define the $C^*$--algebra ${\cal O}_{\rho^g}$. If we use 
$\hat\rho$ to denote the endomorphism of ${\cal O}_{\rho^g}$ induced 
by tensoring on the left by $1_\rho$, then the estimates of Theorem 4.9 
show that $\hat\rho(K)\subset (K^{r+1},K^{r+2})$ in ${\cal O}_{\rho^g}$
and that 
${\cal O}_{\rho^g}$ can be identified canonically with ${\cal O}_K$. Hence 
there is a unitary $R\in(K^{r+2},K^{r+2})$ such that $R\psi=\hat\rho(\psi)$, 
$\psi\in K$. In other words $\hat\rho$ can be identified with $\lambda_R$, 
an algebraic endomorphism of grade zero and rank $\leq r+1$. This estimate 
on the rank is only an upper bound. In fact, if $K$ is ambidextrous, then 
$\hat\rho^g(K)\subset (K,K^2)$. We summarize the discussion as follows.\smallskip

\noindent
{\bf 6.4 Proposition} {\sl Let $\rho$ be an object in a tensor $W^*$--category 
and $K$ a Hilbert space of support one in $(\rho^r,\rho^{r+g})$ with 
$(\rho^r,\rho^r)\times 1_{\rho^{rg}}\subset(K^r,K^r)$. Then 
$({\cal O}_{\rho^g},\hat\rho)= ({\cal O}_K,\lambda_R)$, where 
$R\in(K^{r+2},K^{r+2})$ is the unitary operator in ${\cal O}_\rho$ such 
that $R\psi=\hat\rho(\psi)$, for all $\psi\in K$.}\smallskip

  Returning to the question of when a pair $({\cal A},\hat\rho)$ consisting of a unital $C^*$--algebra and
 a unital endomorphism arises from a system of the form 
$({\cal O}_V, \rho_V)$, where $V$ is a multiplicative unitary, we shall see that the
following conditions are necessary and sufficient: 
\begin{description}
\item{\rm a)} ${\cal A}$ is generated by a tensor $W^*$--category which is 
a tensor subcategory of the category of intertwiners between the powers of 
$\hat\rho$. The generating object will be denoted by $\rho$; 
\item{\rm b)} there is a Hilbert space of support $I$, 
$K\subset(\rho,\rho^2)$;
\item{\rm c)} $\rho(K)\subset(K, K^2)$;
\item{\rm d)} $(\rho,\rho)\subset(K,K)$. \end{description}\smallskip 

We shall see in Theorem 6.11 that ${\cal O}_\rho$ 
is simple, so, in virtue 
of a), $({\cal A},\hat\rho)$ is the system 
$({\cal O}_\rho,\hat\rho)$ associated with the object $\rho$ of the 
tensor $W^*$--category.
\smallskip

  We recognize that c) just says that $\lambda_R$ of Proposition 6.4 
has grade zero and rank $\leq 1$, or equivalently that $R\in(K^2,K^2)$. 
{\rm b)}, on the other hand, tells us that $\lambda_R\lambda_R=\rho_K\lambda_R$, 
or, equivalently that $V:=R\vartheta$ is a multiplicative unitary.
If $V$ is even regular then c) can be strengthened to $K^*\rho(K)=KK^*$.
\smallskip 

  Notice that, as a consequence of Theorem 6.3 and the following 
remark, $({\cal O}_V,\rho_V)$ satisfies a) to d) above relative to the 
Hilbert space $H:=V^*K$.\smallskip   

   As Theorem 6.3 establishes the isomorphism of $({\cal O}_V,\rho_V)$ and 
$({\cal O}_K,\lambda_R)$, $({\cal O}_K,\lambda_R)$ satisfies a) to d),
above, 
too. Indeed, except for d), this is easily seen directly. However, 
$(\lambda_R,\lambda_R)$ being  just the relative commutant of the Hilbert 
space $\lambda_R(K)=RK$, we do  get $(\lambda_R,\lambda_R)\subset(K,K)$ 
by Corollary 5.4 if $V$ is regular. We now prove our first duality 
result.\smallskip

\noindent{\bf 6.5 Theorem} {\sl Let $({\cal A},\hat\rho)$ satisfy 
a) to d) then there is a unique multiplicative unitary $V$ on
the Hilbert space  $K^2\ ,$ such that 
$({\cal A},\hat\rho, K)=({\cal O}_V, \rho_V,V^*K)\ ,$ where 
$V$ is regarded as a representation of $V$. Two systems of the form 
$({\cal O}_V,\rho_V,V^*K)$ for multiplicative unitaries $V$ on $K^2$ are  
isomorphic if and only if the multiplicative unitaries are equivalent.}\smallskip

\noindent{\bf Proof.} By Proposition 6.4, there is a unitary $V\in{\cal A}$ 
such that $({\cal A},\hat\rho,K)=({\cal O}_K,\lambda_{V\vartheta_{K,
K}},K)$. 
But c) implies that $V\in (K^2,K^2)$ and b) may be read as 
$\hat\rho^2=\rho_K\hat\rho$ so that $V$ is a multiplicative unitary on 
$K^2$ \cite{Cu2},\cite{R}. By Theorem 6.3 again, we can consider
isomorphisms 
of systems of the form $({\cal O}_K,\lambda_R,K)$. So if $\tau$ is 
an isomorphism from ${\cal O}_K$ to ${\cal O}_{K'}$ with $\tau(K)=K'$ and 
$\tau\circ\lambda_R=\lambda_{R'}\circ\tau$. Then $\tau(R)=R'$ and, since 
$\tau(\vartheta_{K,K})=\vartheta_{K',K'}$, 
$\tau(V)=V'$. So if $U$ is the unitary from $K$ to $K'$ such that $\tau(\psi)=U\psi$, 
for $\psi\in K$, its second tensor power will intertwine $V$ and $V'$, 
realizing the desired equivalence. The converse is obvious.  
\smallskip

  In the case of a regular multiplicative unitary we can even obtain a 
categorical rather than an algebraic duality theorem. The following 
result characterizes tensor $W^*$--categories of the form ${\cal T}_V$ 
for a regular multiplicative unitary $V$.\smallskip 

\noindent
{\bf 6.6 Theorem} {\sl Given a tensor $W^*$--category ${\cal T}_\rho$
whose objects 
are the tensor powers of a $C^*$-amenable object $\rho$, suppose there is 
an ambidextrous Hilbert space $K$ of support one in $(\rho,\rho^2)$ 
with $(\rho,\rho)\times 1_\rho\subset(K,K)$ and 
$$K^*\times 1_\rho\circ 1_\rho\times K=K\circ K^*$$ 
then the unitary $V$ on $K^2$ defined by 
$$V\circ\psi_1\times 1_\rho\circ\psi_2=1_\rho\times\psi_2\circ\psi_1,\quad 
\psi_1,\,\psi_2\in K,$$ 
is multiplicative and ${\cal T}_\rho$ and ${\cal T}_V$ are
isomorphic.}\smallskip 

\noindent
{\bf Proof.} The fact that $V$ is multiplicative follows from Theorem 4.2 and 
the condition $K^*\times 1_\rho\circ 1_\rho\times K=K\circ K^*$ just 
says that $K$ is regular. We now consider the image of ${\cal T}_\rho$ 
in ${\cal O}_\rho$ under the canonical map, then the conditions a) to d) 
above are satisfied by ${\cal O}_\rho$ and we conclude from Theorem 6.5
that 
$({\cal O}_\rho,\hat\rho,K)=({\cal O}_V,\rho_V,V^*K)$. Now $\rho$ is 
$C^*$--amenable by assumption and $V$ is $C^*$--amenable by Theorem 6.2
and, with
 an obvious notation, ${\cal T}_{\hat\rho}$ and ${\cal
T}_{\rho_V}$ 
are isomorphic hence ${\cal T}_\rho$ and ${\cal T}_V$ are
isomorphic.\smallskip 
  
  We now want to give an algebraic characterization of the situation where
the multiplicative unitary $V$ on $K^2$ is endowed
with a standard braided symmetry $\varepsilon$, a concept explained in the 
appendix. In this case the system 
$({\cal O}_{\hat V}, \rho_{\hat V}, {\hat V}^*K)\ ,$ with 
$\varepsilon=V\vartheta_{K, K}\hat V\ ,$
satisfies conditions a) to d). As pointed out in the appendix, $({\cal O}_{\hat V},
\rho_{\hat V})=({\cal O}_V, \rho_V)\ ,$ where  $V$ is the regular
corepresentation and $\hat V$ the regular representation. A direct computation
shows that the condition for being standard, namely that $\hat V$ and $V_{23}$
commute on $(K^3, K^3)\ ,$ is equivalent to
$\rho_V(\varepsilon\psi)=\rho_H(\varepsilon\psi)$, $\psi\in H:={\hat V}^*K $. 
We start with the following result.\smallskip

 \noindent{\bf 6.7 Theorem} {\sl Let $({\cal A},\rho)$ be a pair 
consisting of a $C^*$--algebra and a unital endomorphism satisfying 
a) to d). Let $\varepsilon$ be a braided symmetry for $\rho$ satisfying 
\begin{description} 
\item{\rm e)} $\rho(\varepsilon\psi)=\rho_K(\varepsilon\psi)\ ,\quad\psi\in K\ .$
\end{description} Then there is a multiplicative unitary $V$ on a Hilbert space
$H$ whose regular corepresentation $V$ is endowed with a standard braided 
symmetry $\varepsilon_V$ of $V$ and an
isomorphism $\Phi: {\cal A}\to{\cal O}_V$ such that 
\begin{description}
\item{\rm 1)} $\Phi\circ\rho=\rho_V\circ\Phi\ ;$
\item{\rm 2)} $\Phi(\varepsilon)=\varepsilon_V\ ;$
\item{\rm 3)} $\Phi(K)=\hat V^*H\ ,$ where $\varepsilon_V=V\vartheta\hat V\ .$
\end{description}
The multiplicative unitary $V$ is
determined, up to equivalence, by  the above conditions.
If $\varepsilon$ is a permutation symmetry then $({\cal A}, \rho)$
correponds to a locally compact group $G$ and $\varepsilon$ to the usual
permutation symmetry. }\smallskip

\noindent{\bf Proof.} By Theorem 6.5 we can find a regular
representation $\hat V$ acting on a Hilbert space $H$ such that the triples
 $({\cal
A}, \rho, K)$ and $({\cal O}_{\hat V}, \rho_{\hat V}, {\hat V}^*H)$ are isomorphic
via an isomorphism $\Phi\ .$ If we write
$\varepsilon_V:=\Phi(\varepsilon)=V\vartheta\hat V$ then by e) $V_{23}$ and
$\hat V$ commute on $H^3\ ,$ so as shown in the appendix, $V$ is multiplicative and
$\varepsilon_V$ is a standard braided symmetry of $V\ .$ It will follow 
from the next proposition that $V$ is unique up to equivalence. If
$\varepsilon$ is a permutation symmetry then, by Proposition A.4, $V$ is
cocommutative, thus coming from a locally compact group $G$ \cite{BS}. 
\smallskip

\noindent{\bf 6.8 Proposition.} {\sl Let $V$ and $V'$ be multiplicative unitaries on
Hilbert spaces $H^2$ and ${H'}^2$  endowed with standard braided symmetries
$\varepsilon_V=V\vartheta_{H, H}\hat V$ and
$\varepsilon_{V'}=V'\vartheta_{H', H'}\hat{V'}$ respectively. If there is an
isomorphism $\Phi: {\cal O}_V\to{\cal O}_{V'}$ satisfying \begin{description}
\item{\rm a)} $\Phi\circ\rho_V=\rho_{V'}\circ\Phi\ ;$
\item{\rm b)} $\Phi(\varepsilon_V)=\varepsilon_{V'}\ ;$
\item{\rm c)} $\Phi(K)=K'\ ,$ where $K=\hat V^*H$ and $K'= \hat{V'}^*H'\ ,$
\end{description}
then $V$ and $V'$ are unitarily equivalent.}\smallskip

\noindent{\bf Proof.} By c), $\Phi\circ\rho_K=\rho_{K'}\circ\Phi$ on ${\cal O}_V$.
Thus by a) and b) $\Phi\circ\lambda_{\hat V^*}(V)=\lambda_{\hat {V'}^*}(V')$ since
$\lambda_{\hat V^*}(V)=\rho_V(\varepsilon^*_V)\rho_K(\varepsilon_V)$ with a
similar result for $V'\ .$ Now by c) there is a unitary operator $U: H\to H'$ extending to
an isomorphism $\alpha: {\cal O}_H\to{\cal O}_{H'}$ such that 
$\Phi\circ\lambda_{\hat V^*}=\lambda_{\hat {V'}^*}\circ\alpha$, on ${\cal O}_H$. It follows that
$\alpha(\hat V)=\hat{V'}\ .$ Let us define $\tilde{\varepsilon}_V=\hat
V\varepsilon_V\hat V^*=\hat VV\vartheta_{H, H}\ .$ Then
$\varepsilon_V=\lambda_{\hat V^*}(\tilde{\varepsilon}_V)$ since $(V^{\times 2},
V^{\times 2})$ is generated by $K(V, V)K^*$ as a weakly closed subspace and $\hat
V$ commutes with $(V, V)\ .$ If we define the operator
$\tilde{\varepsilon}_{V'}$ corresponding to $V'\ ,$ as above, then  we deduce 
$\alpha(\hat VV\vartheta_{H, H})=\hat{V'}V'\vartheta_{H', H'}\ ,$ from b) and the
previous  relation intertwining $\Phi$ and $\alpha$. 
Clearly $\alpha(\vartheta_{H, H})=\vartheta_{H', H'}$, so $\alpha(V)=V'$. Now 
the adjoint action of $U\times U$ on $(H^2, H^2)$  implements $\alpha$, 
completing the proof. \smallskip

\noindent
{\bf Remark.} We can now complement the discussion following Theorem 6.2. 
With a standard braided symmetry, ${\cal C}(V)$ and ${\cal R}(V)$ are 
isomorphic as tensor $W^*$--categories embedded in Hilbert spaces. 
Hence by Theorem 6.2, a  $V$ with a standard braided symmetry and 
$\hat V$ regular is $C^*$--amenable in ${\cal C}(V)$. 
Indeed $\hat V^*K\subset (V,V^2)$ has trivial relative commutant in ${\cal
O}_K$.
\smallskip

   We now examine model endomorphisms which are more $C^*$--algebraic in nature.
In fact, when $K$ is infinite dimensional, the above systems are not really  $C^*$--algebraic 
in nature since $\O_K$ is the norm closure of subspaces endowed with a 
$W^*$--topology and is not even separable when $K$ is an infinite dimensional 
separable Hilbert space. To cure this defect, we might replace 
$({\cal O}_K,\lambda_R)$ by $({\cal P}_K,\tau_R)$, where ${\cal P}_K$ is
the 
smallest $C^*$--subalgebra of ${\cal O}_K$ containing $K$ and stable under 
$\lambda_R$ and $\tau_R$ denotes the restriction of $\lambda_R$ to ${\cal
P}_K$. 
Note that, since $K\subset{\cal P}_K$, 
$$(\tau_R^m,\tau_R^n)=(\lambda_R^m,\lambda_R^n)\cap{\cal P}_K.$$ 
At the same time, we would now like our model endomorphisms to have the form 
$({\cal O}_\rho,\hat\rho)$ where $\rho$ is an object in a tensor
$C^*$-category. 
We still want a Hilbert space $K$ of intertwiners between powers of $\rho$. 
We say that a Hilbert space $H\subset (\rho,\sigma)$ in a $C^*$--category has 
{\it zero left annihilator} if $X\circ\psi=0$ for all $\psi\in H$ implies
$X=0$. It obviously suffices to take $X$ to be a positive element of 
$(\sigma,\sigma)$.\smallskip 

  We first prove a result that will imply that the $C^*$--algebras of interest 
are simple $C^*$--algebras.\smallskip 

\noindent
{\bf 6.9 Theorem} {\sl Let ${\cal T}_\rho$ be a tensor $C^*$--category 
whose objects are the powers of the object $\rho$. Let $K\in(\rho^r,\rho^{r+g})$ 
be a Hilbert space such that $K\times 1_{\rho^m}$ has left annihilator 
zero for $m\in{\Bbb N}_0$ and suppose that 
$$(\rho^r,\rho^r)\times 1_{\rho^r}\subset(K^r,K^r).$$ 
Then ${\cal O}_\rho$ is a simple $C^*$--algebra.}\smallskip 

\noindent
{\bf Proof.} The proof follows that of the simplicity of the extended 
Cuntz algebra Theorem 3.1 of \cite{CDPR}. Obviously, $K$ must now 
play the role of the generating Hilbert space. We identify the arrows of 
${\cal T}_\rho$ with their images in ${\cal O}_\rho$. It suffices to 
show that any non-degenerate representation $\pi$ of ${\cal O}_\rho$ is 
faithful. $\pi$ is trivially isometric on Hilbert spaces in ${\cal
O}_\rho$ 
and in particular on $K^n$, $n\in{\Bbb N}_0$. But then $\pi$ is also 
isometric on $(K^m,K^n)$, defined as in the discussion preceding
Lemma~4.11. 
Given $k\in{\Bbb Z}$, we let 
$${}^o{{\cal O}_\rho}^k:=\cup_k(K^r,K^{r+k}),\quad r\geq 0,\,\,r+k\geq 0$$ 
and let ${}^o{\cal O}_\rho$ denote the $^*$--subalgebra of ${\cal O}_\rho$ 
obtained by taking finite sums of elements from the 
${}^o{{\cal O}_\rho}^k$. We 
have seen that $\pi$ is isometric on each ${}^o{\cal O}^k_\rho$. If 
$K$ is infinite dimensional, we can continue as in the proof of Theorem~3.1 
of \cite{CDPR} to show that $\pi$ is isometric on ${}^o{\cal
O}_\rho$. However 
this algebra is dense in ${\cal O}_\rho$, since arguing as in Theorem 4.9, 
we have $(\rho^{r+m},\rho^{r+n})\times 1_{\rho^r}\subset(K^{r+m},K^{r+n})$. 
Hence $\pi$ is isometric and ${\cal O}_\rho$ is simple. If $K$ is finite 
dimensional, we need only remark that $(K^m,K^n)$ is the set of all 
linear mappings from $K^m$ to $K^n$, so that, as a $C^*$--algebra, 
${\cal O}_\rho={\cal O}_K$ which is a simple $C^*$--algebra.

  Thus instead of considering $(\O_V,\rho_V)$, we consider the smallest 
$C^*$--subalgebra $\A_V$ of $\O_V$ containing $H:=V^*K$ and stable under 
$\rho_V$, equipped with the endomorphism $\sigma_V$ obtained by restricting  
$\rho_V$. The system $(\A_V, \sigma_V)$ is then a natural 
candidate for such a minimal $C^*$--model system and we therefore address 
the question of finding 
necessary and sufficient conditions on a system $(\A, \hat\sigma)$ for it to be 
of the form $(\A_V, \sigma_V)$. We shall only discuss the case where $V$ is 
its regular representation, as the case of a corepresentation 
can, as before, be reduced to this case, for systems endowed with 
a standard braided symmetry. (This is not completely trivial, in that
the concept of 
a  braided symmetry not taking values in the intertwiner spaces, but just in
their weak closures in some Hilbert space representation with support $I$ 
has to be formalized.)

We start by pointing out that the smallest tensor $C^*$--subcategory
${\cal S}_V$ 
of $\cT_V$ containing $V$ and the intertwining space $H=V^*\circ K\times 1_K$ is, up to tensoring 
on the right by $1_V$, $W^*$--dense in ${\cal T}_V$. 
More precisely, $(V^r, V^s)\times1_V$ is contained in the $W^*$--closure of $H^s{H^r}^*$ in 
$(V^{r+1}, V^{s+1})$. Note that $(\A_V, \sigma_V)$ is the canonical system 
derived from $\sigma=V$ regarded as an object of the tensor $C^*$--category 
${\cal S}_V$. It is therefore canonically associated
with $\cT_V$. The following simple result
relates the systems associated with $V$ 
considered as an object of  ${\cal S}_V$ and ${\cal T}_V$,
respectively. \smallskip

\noindent{\bf 6.10 Proposition} {\sl For $r,s=0,1,2\dots$, 
$({\sigma_V}^r,{\sigma_V}^s)=({\rho_V}^r, {\rho_V}^s)\cap\A_V$}\smallskip

\noindent{\bf Proof.} As $H$ generates ${\cal O}_V$, $T\in({\rho_V}^r,
{\rho_V}^s)$ if it
 intertwines the restrictions of the corresponding endomorphisms to the space
of intertwiners $H=V^*K$. Now this space is contained in $\A_V$,
therefore $({\rho_V}^r, {\rho_V}^s)\cap\A_V \supset({\sigma_V}^r,{\sigma_V}^s)$, 
and the reverse inclusion is obvious.\smallskip

Summarizing, the above  discussion and conditions a)--d) lead to  
the following  necessary conditions for
$(\A, \hat{\sigma})$ to be of the form $(\A_V, \sigma_V)$ with $V$ a {\it regular} 
multiplicative unitary:
\begin{description}
\item{\rm a$'$)} 
${\cal A}$ is the smallest $\hat\sigma$--stable $C^*$--subalgebra
containing a 
Hilbert space $K$ with zero left annihilator,
\item{\rm b$'$)} $K\subset(\hat\sigma,\hat\sigma^2)$,
\item{\rm c$'$)} $K^*\hat\sigma(K)=KK^*$,
\item{\rm d$'$)} $(\hat\sigma,\hat\sigma)K=K$. 
\end{description} 

  To see that these conditions are necessary, we need only remark that 
c$'$) just expresses the regularity of the multiplicative unitary and 
ensures that the hypotheses of Lemma 5.3 hold for $L=M=\hat\sigma(K)$ as we 
have already remarked.  d$'$) is now a consequence of Corollary 5.4. 
The sufficiency of these conditions follows from the next result.\smallskip

\noindent{\bf 6.11 Theorem} {\sl Let $(\A, \hat\sigma, K)$ satisfy conditions a$'$) to d$'$)
then ${\cal A}$ is simple and there is a regular multiplicative unitary
$V$, unique up to equivalence,
 such that $(\A,\hat\sigma,K)$ is isomorphic to the model  dual object
$(\A_V,  \sigma_V, V^*K)$, where $V$ is regarded as a representation of $V$.}\smallskip

\noindent{\bf Proof.}  ${\cal A}$ is simple by Theorem 6.9. From condition
c$'$) 
it follows that $$K^*K^*\hat\sigma(K)K\subset{\Bbb C}$$ 
and hence, since $K$ has left annihilator zero that 
$\hat\sigma(K)\subset (K,K^2)$. We now may apply Lemma~4.11 and
Theorem~4.12 to the 
tensor $C^*$--category of intertwiners between the powers of $\hat\sigma$ 
to conclude that there is a multiplicative unitary $V$ on $K^2$ with 
$V\theta\psi=\hat\sigma(\psi)$, $\psi\in K$. 
$V$ is regular by c$'$), therefore the remaining
conclusions follow using arguments similar to those of Theorem~6.5.
\smallskip 

  We can also give necessary and sufficient conditions for $({\cal 
A},\hat\sigma)$ 
to be of the form $({\cal A}_V,\sigma_V)$ for a general multiplicative
unitary. 
We retain a$'$) and b$'$) above and replace c$'$) and d$'$) by 
\begin{description}
\item{c$''$)} $\hat\sigma(K)\subset(K,K^2)$,
\item{d$''$)} If we consider the tensor $C^*$--subcategory ${\cal T}$ of
the tensor 
$C^*$--category of intertwiners between the powers of $\hat\sigma$ generated 
by $K$ and denote its objects by $\sigma^n$, where $n\in{\Bbb N}_0$, then 
$(\sigma, \sigma)K=K$.
 
\end{description}

\noindent{\bf 6.12 Theorem} {\sl Let $(\A, \hat\sigma, K)$ satisfy conditions 
a$'$), b$'$), c$''$) and d$''$)
then ${\cal A}$ is simple and there is a multiplicative unitary $V$,
unique up to equivalence,
 such that $(\A,\hat\sigma,K)$ is isomorphic to the model  dual object
$(\A_V,  \sigma_V, V^*K)$, where $V$ is regarded as a representation of $V$.}\smallskip

   The proof is a simplified version of that of the previous theorem seeing 
that regularity now plays no role. Finally, as a pendant to Theorem~6.6, we 
give a characterization of tensor $C^*$--categories of the form ${\cal 
S}_V$ 
for a multiplicative unitary $V$\smallskip 

\noindent
{\bf 6.13 Theorem} {\sl Let ${\cal T}_\rho$ be a tensor $C^*$--category
whose 
objects are the tensor powers of an object $\rho$ and suppose ${\cal 
T}_\rho$ 
is generated by an ambidextrous Hilbert space $K$ in $(\rho,\rho^2)$ 
such that $K\times 1_{\rho^m}$ has left annihilator zero for $m\in{\Bbb
N}_0$. 
Let $V$ be the associated multiplicative unitary of Theorem 4.12, then 
${\cal T}_\rho,\,K$ is isomorphic to ${\cal S}_V,\,V^*\circ K\times
1_\rho$.}\smallskip 

\noindent
{\bf Proof.} We let $\tilde{\cal K}$  be the category of Hilbert spaces 
with commuting shifts $F$ and $G$ associated with the ambidextrous space 
$K\in(\rho,\rho^2)$ as in Theorem 4.12. Then the functor $G_{V^*}$ 
not only commutes with $F$ but satisfies $G_{V^*}G=\hat FG_{V^*}$ by 
Proposition 3.1e.  If we interpret $G_{V^*}$ as a tensor $^*$--functor 
as in Proposition 4.10, we recognize that the image 
of ${\cal T}_\rho$ is the tensor $C^*$--subcategory generated by 
$H:=V^*\circ K\times 1_\rho\in{\cal T}_V$. But this is, by definition, 
the minimal $C^*$--subcategory ${\cal S}_V$.

\section{Appendix. Braided Symmetry}

  This appendix is devoted to the notion of braided symmetry. This 
evolved from a notion of the same name introduced in \cite{CDPR} in the 
context of endomorphisms of $C^*$--algebras to generalize a previous more 
restricted notion of permutation symmetry in \cite{DR}. Expressed in the context 
of a general tensor category, this notion can be expressed as follows. 
Let ${\sigma}$ denote the
endomorphism of the braid group ${\mathbb B}_{\infty}=\cup{\mathbb B}_n$
that
shifts
the braids $b\in{\mathbb B}_n$ on $n$ threads to the right.  By a braided
symmetry for an object  $V$ in a tensor category ${\cal T}$  we mean a 
representation $\varepsilon$ of
${\mathbb B}_{\infty}$ in ${\cal T}$ such that 
\begin{eqnarray*} 
&&\varepsilon(b)\in(V^{\times
n}, V^{\times n})\ ,\quad b\in{\mathbb B}_n\ ,\\ 
&&\varepsilon(\sigma(b))=1_V{\times}\varepsilon(b)\
,\quad b\in{\mathbb B}_\infty=\cup{\mathbb B}_n\ ,\\ 
 &&\varepsilon(s,1)\circ X\times 1_V=1_V\times X\circ  \varepsilon(r,1)\ ,\quad
X\in(V^{\times r}, V^{\times s})\ ,\end{eqnarray*} 
where  $(1, 1)=b_1$ is the braid on the first two
threads and $(s, 1)=b_1\sigma(b_1)\dots\linebreak\sigma^{s-1}(b_1)\,.$

Obviously, if the full subcategory whose objects are the tensor powers of V can 
be made into a braided tensor category then this braiding does define a braided 
symmetry. However, not all braided symmetries arise in this way. An 
application of this notion of braided symmetry is Theorem 5.31 of 
\cite{LR} relating notions of amenability to the existence of a unitary 
braided symmetry.\smallskip 

  Whilst this definition, with a view to simplicity, focused attention on 
the full subcategory whose objects are the tensor powers of $V$, we here need 
to consider the full tensor category and we will hence call ${\cal T}$ braided 
relative to a distinguished object $V\in{\cal T}$ if for any object $W$ in 
${\cal T}$ there is an invertible arrow $\var_W\in(W\times V, V\times W)$ such that
\begin{eqnarray*}
&&\var_{W\times W'}=\var_W\times 1_{W'}\circ1_W\times\var_{W'}\ ,\\
&&\var_{W'}\circ T\times 1_V=1_V\times T\circ\var_W\ ,\quad T\in(W, W')\ .
\end{eqnarray*}
The second equation just says that $\var$ is a natural transformation from 
the functor of tensoring on the right by $V$ to that of tensoring on the left by $V$. 
The first equation implies in particular that this natural transformation takes the 
value $1_V$ on the tensor unit. If $\var$ is a braided symmetry for $V$ then 
we may define braided symmetries for the tensor powers $V^{\times n}$ of $V$ 
inductively, setting 
$$\var_W^{\times n}:= 1_V\times\var_W^{\times n-1}\circ\var_W\times 1_{V^{n-1}}.$$ 
It is easy to see that  if $\var_W$ is defined on a 
subset of objects, closed under tensor products, and such that every object of 
${\cal T}$ is a subobject of a (finite) direct sum of objects from the
subset and the equations 
are satisfied for $W$ and  $W'$ in the subset, then $\var_W$ extends uniquely to a 
braided  symmetry on the whole category. This remains true if ${\cal T}$
is a 
$W^*$--category and we allow infinite direct sums.\smallskip 

   Braided symmetries in this sense have appeared 
as the starting point of the centre construction in tensor categories, 
see eg.\ \cite{K} where references to the original articles are given. 
But in view of its simplicity, the notion may well have appeared in 
other contexts, still unknown to the authors. We recall, however, 
the notion of an arrow between braided symmetries. If $\var$ and $\var'$ are
braided symmetries for $V$ and $V'$, respectively, then an arrow 
$T\in(\var,\var')$ is an arrow $T\in(V,V')$ such that 
$$\var'_W\circ 1_W\times T=T\times 1_W\circ\var_W,$$ 
for each object $W$ of ${\cal T}$.\smallskip
 
When does an arrow $\var_V\in(V\times V,V\times V)$ define a 
braided symmetry? To give some kind of answer, we consider a strict tensor 
category ${\cal T}$ and an object $V$ with the property that the functor
of 
tensoring on the right by $V$ is faithful and such that $W\times V$ is a
direct sum of copies of $V$ for each object $W$ of ${\cal T}$. This
property 
is related to the notion {\it right regular} representation. The 
corresponding property of $V\times W$ being a direct sum of copies of $V$ 
is similarly related to the notion of  
{\it left regular} representation.\smallskip

Given an invertible $\var_V\in (V\times V,V\times V)$ such that 
$$\var_V\circ T\times 1_V=1_V\times T\circ\var_V,\quad T\in(V,V),$$
there is, for each object $W$ of ${\cal T}$ a unique invertible 
$$\var_{W\times V}\in (W\times V\times V,V\times W\times V)$$ 
such that 
$$\var_{W\times V}\circ S\times 1_V=1_V\times S\circ\var_V,\quad S\in(V,W\times V).$$
In fact, we have only to pick $X_i\in(V,W\times V)$ and $Y_i\in(W\times V,V)$ 
such that $\sum_iX_iY_i=1_{W\times V}$ and we see that we have no option but
to set 
$$\var_{W\times V}:=\sum_i1_V\times X_i\circ\var_V\circ Y_i\times 1_V.$$ 
A routine computation shows that
$$\var_{W'\times V}\circ S\times 1_V=1_V\times S\circ \var_{W\times V},\quad
S\in (W\times V,W'\times V).$$ 
Now consider the set $\Sigma$ of objects $W$ such that
$$\var_{W\times V}\circ 1_W\times\var_V^{-1}\in (W\times V,V\times W)\times 1_V$$
and define $\var_W$ by
$$\var_W\times 1_V=\var_{W\times V}\circ 1_W\times\var_V^{-1}.$$
If $W$ and $W'$ are in $\Sigma$ and $T\in (W,W')$,
$$\var_{W'}\times 1_V\circ T\times 1_{V\times V}=\var_{W'\times V}\circ 1_{W'}
\times\var_V^{-1}\circ T\times 1_{V\times V}$$
$$=\var_{W'\times V}\circ T\times 1_{V\times V}\circ 1_W\times\var_V^{-1}
=1_V\times T\times 1_V\circ\var_{W\times V}\circ 1_W\times\var_V^{-1},$$
so that
$$\var_{W'}\circ T\times 1_V=1_V\times T\circ\var_W,\quad T\in(W,W').$$
The set $\Sigma$ trivially contains the tensor unit. Suppose both $W$ and $W'$ 
are in $\Sigma$ then 
$$\var_{W\times W'\times V}\circ 1_W\times\var_{W'\times V}^{-1}\circ 
1_{W\times V}\times S=\var_{W\times W'\times V}\circ 1_W\times S\times 
1_V\circ 1_W\times\var_V^{-1}$$
$$=1_V\times 1_W\times S\circ\var_{W\times V}\circ 1_W\times \var_V^{-1}$$ 
$$=1_V\times 1_W\times S\circ \var_W\times 1_V=\var_W\times 1_{W'\times V}\circ 
1_W\times 1_V\times S,$$ 
where $S\in (V,W'\times V)$. Now since $W'\times V$ is 
a direct sum of copies of $V$ we conclude that
$$\var_{W\times W'\times V}\circ 1_W\times\var_{W'\times V}^{-1}=\var_W\times 1_{W'\times V}.$$ 
But since $W'\in\Sigma$, $\var_{W'\times V}=\var_{W'}\times 1_V\circ 1_{W'}\times\var_V$. 
Thus 
$$\var_{W\times W'\times V}=\var_W\times 1_{W'\times V}\circ 1_W\times\var_{W'}\times 1_V,$$
so that $W\times W'\in\Sigma$. Furthermore, 
$$\var_{W\times W'}=\var_W\times 1_{W'}\circ 1_W\times\var_{W'}.$$ 
It is also easy to see that $\Sigma$ is closed 
under subobjects and direct sums.

We have now proved the following result.\smallskip 

\noindent{\bf A.1 Proposition} {\sl Let ${\cal T}$ be a tensor category,
$V$ an 
object such that $W\times V$ is a direct sum of copies of $V$ for each $W$. 
Let $\var_V\in (V\times V,V\times V)$ be a unitary such that 
$$\var_V\circ T\times 1_V= 1_V\times T\circ\var_V,\quad T\in(V,V),$$
$$\var_V\times 1_V\circ 1_V\times\var_V\circ S\times 1_V=1_V\times S\circ \var_V,
\quad S\in(V,V\times V).$$
then there is a unique maximal tensor subcategory with a braided 
symmetry $\var$ whose value at $V$ coincides with the given invertible $\var_V$.
This is a full subcategory closed under subobjects and direct sums.}\smallskip

\noindent{\bf Proof.} We need only remark that the second equation above implies 
that 
$$\var_{V\times V}= \var_V\times 1_V\circ 1_V\times\var_V$$
and hence $V$ is in $\Sigma$ and the two definitions of $\var_V$ coincide.\smallskip

   We now consider the category ${\cal C}(V)$ of corepresentations 
of a multiplicative unitary $V$ with its forgetful functor $\iota$ into the  
underlying category of Hilbert spaces. We let $\vartheta$ denote the 
braided symmetry relative to $\iota(V)$ derived from the symmetry on the 
category of Hilbert spaces but write $\vartheta_W$ in place of $\vartheta_{\iota(W)}$. 
The relevance of braided symmetries to this paper lies in the fact that 
there are many cases where ${\cal C}(V)$ admits a braided symmetry $\var$ 
relative to $V$ with the further property that $\hat V$ defined by 
$V\vartheta\hat V=\var_V$ is another multiplicative unitary on the same space. Such a braided symmetry will be called 
{\it standard}. $\hat V$ determines and is uniquely determined by $\var$.\smallskip 

\noindent
{\bf Theorem A.2} {\sl A braided symmetry $\var$ on ${\cal C}(V)$ relative
to 
$V$ is standard if and only if 
$$\hat V_{12}V_{23}=V_{23}\hat V_{12}.$$ 
If $\var$ is standard and $W$ is a corepresentation of $V$, then $\hat W$ 
defined by 
$$W\vartheta_{W,K}\hat W=\var_W$$ 
is a representation of $\hat V$ and for any pair $W$, $W'$ of corepresentations, 
we have 
$$\hat W_{12}W'_{23}=W'_{23}\hat W_{12},$$ 
$$\widehat{W\times W'}=\hat W\times\hat W'.$$}

\noindent
{\bf Proof.} Whether $\var$ is standard or not, a simple calculation shows 
that $T\in(W,W')$ if and only if $T\in(\hat W,\hat W')$. Since 
$\hat W\in(W\times V,\iota(W)\times V)$, this yields 
$$\hat W_{12}\widehat{W\times V}=\widehat{\iota(W)\times V}\hat W_{12}.$$ 
Now 
$$\hat W_{13}\hat W'_{23}=\hat W_{13}\vartheta_{23}W^{'-1}_{23}\var_{23}$$
$$=\vartheta_{23}\hat W_{12}W^{'-1}_{23}\var_{23}.$$ 
If $\hat W_{12}W'_{23}=W'_{23}\hat W_{12}$ then we get 
$$\hat W_{13}\hat W'_{23}=\vartheta_{23}W^{'-1}_{23}\vartheta^{-1}_{12}W^{-1}_{12}
\var_{W\times W'}=\vartheta^{-1}_{W\times W'}W^{'-1}_{13}W^{-1}_{12}\var_{W\times W'}.$$ 
So $\hat W_{13}\hat W'_{23}=\widehat{W\times W'}$. Thus if $\hat W_{12}$ 
and $V_{23}$ commute, we could conclude from our first identity that 
$$\hat W_{12}\hat W_{13}\hat V_{23}=\hat V_{23}\hat W_{12},$$ 
and, in particular, if $\hat V_{12}$ and $V_{23}$ commute that 
$\hat V$ is a multiplicative unitary.
We now show that $\hat V_{12}V_{23}=V_{23}\hat V_{12}$ implies 
$\hat W_{12}W'_{23}=W'_{23}\hat W_{12}$. Now
$$V_{23}\hat V_{12}W_{24}W_{34}=\hat V_{12}V_{23}W_{24}W_{34}=
\hat V_{12}W_{34}V_{23},$$
$$V_{23}W_{24}\hat V_{12}W_{34}=V_{23}W_{24}W_{34}\hat V_{12}=
W_{34}V_{23}\hat V_{12},$$  
where we have used the multiplicativity of $V$.
But these expressions are equal, so cancelling, we find 
$\hat V_{12}W_{24}=W_{24}\hat V_{12}$ and we have replaced $V$ by 
$W$. The proof is completed by a similar step using the multiplicativity of $\hat V$.
$$\hat W_{12}\hat W_{13}W'_{34}\hat V_{23}=\hat W_{12}\hat W_{13}\hat 
V_{23}W'_{34}=\hat V_{23}\hat W_{12}W'_{34},$$
$$\hat W_{12}W'_{34}\hat W_{13}\hat V_{23}=W'_{34}\hat W_{12}
\hat W_{13}\hat V_{23}=W'_{34}\hat V_{23}\hat W_{12}.$$ 
Again these two expressions are equal and cancelling gives 
$\hat W_{13}W'_{34}=W'_{34}\hat W_{13}$, as claimed. Hence, our previous 
computation shows that $\hat W$ is a representation of $\hat V$, 
completing the proof.\smallskip  

\noindent
{\bf Remark.} Since a braided symmetry relative to $V$ might not be defined 
on the whole of ${\cal C}(V)$, it is worth remarking that ${\cal C}(V)$ 
can be replaced by a full tensor subcategory in the above theorem.\smallskip

  Note that the above theorem shows that when $\var$ is a standard 
braided symmetry, ${\cal C}(V)$ and ${\cal R}(\hat V)$ are canonically 
isomorphic as tensor $W^*$--categories. Of course, we could 
start with the multiplicative unitary $\hat V$ and then use $\var$ to 
define $V$ and $V$ would again be multiplicative if and only if 
$\hat V_{12}V_{23}=V_{23}\hat V_{12}$. 
We conclude by showing that the 
interesting standard braided symmetries cannot be permutation symmetries.\smallskip

\noindent{\bf A.3 Lemma} {\sl Let $\varepsilon=V\vartheta\hat V$ define a
standard braided symmetry of $V\ .$ Then $\hat V^{-1}\hat V^{-1}_{23}
\hat V=\varepsilon^{-1}_{23}{\hat V}^{-1}\varepsilon_{23}\hat V^{-1}_{23}\ .$}\smallskip

\noindent{\bf Proof.} 
$$\varepsilon^{-1}_{23}{\hat V}^{-1}\varepsilon_{23}^{-1}\hat
V_{23} = (V\vartheta\hat V)^{-1}_{23}{\hat V}^{-1}(V\vartheta)_{23} =$$
$$\hat V^{-1}_{23}\hat V^{-1}_{13} = \hat V^{-1}_{12}\hat
V^{-1}_{23}\hat V_{12} \ ,$$
since $V_{23}$ and $\hat V$ commute and $\hat V$ is multiplicative.
\smallskip

\noindent{\bf A.4 Proposition} {\sl Let $\varepsilon=V\vartheta\hat V$ be a
standard braided symmetry for $V\ .$ Then $\vartheta=\varepsilon
\hat V^{-1}_{23}\varepsilon_{23}\hat
V(\varepsilon^{-1}_{23})^2\hat V^{-1}\varepsilon_{23}\hat V_{23}\ .$ In particular if
$\varepsilon$ is a permutation symmetry then
$\varepsilon=\vartheta$ and $V$ is cocommutative.}\smallskip

\noindent{\bf Proof.}
 $(V^{\times 2}, V^{\times 2})$ is generated, as a weakly closed subspace,
by $H(V, V)H^*\ ,$ with $H=\hat V^{-1}K\ ,$ thus setting $\tilde{\varepsilon}=\hat
V\varepsilon\hat V^{-1}=\hat V V\vartheta$ then $\varepsilon=\hat V^{-1}
\hat V^{-1}_{23}\tilde\varepsilon\ \hat V_{23}\hat V.$  So $\varepsilon$ 
commutes with $\hat V^{-1}\hat V^{-1}_{23}\hat V$. Thus by the previous lemma
$$\varepsilon =\hat V_{23}\varepsilon^{-1}_{23}{\hat V}\varepsilon_{23}
\varepsilon\varepsilon^{-1}_{23}{\hat V}^{-1}\varepsilon_{23}\hat V^{-1}_{23}\ .$$ 
Now $\varepsilon$ defines a braided symmetry for
$V$, thus for  any $\pp=\hat V^{-1}\ff\in H\subset(V, V^{\times 2})\ ,$
$\varepsilon\varepsilon_{23}\pp\varepsilon^{-1} = (\hat V^{-1}\ff)_{23}=
\hat V^{-1}\hat V^{-1}_{23}\hat V\hat V_{13}\vartheta\ff$  since $\hat V$ is multiplicative, and,
again by Lemma A.3, $ \varepsilon^{-1}_{23}{\hat V}^{-1}\varepsilon_{23}\hat V^{-1}_{23}
\vartheta\hat V_{23}=\varepsilon\varepsilon_{23}\hat V^{-1}\varepsilon^{-1}_{23}$
 hence 
$$\vartheta=\hat V_{23}\varepsilon^{-1}_{23}\hat V\varepsilon_{23} 
\varepsilon\varepsilon_{23}\hat V^{-1}\varepsilon^{-1}_{23}\hat V^{-1}_{23}$$
and the conclusion follows.
\end{document}